\documentclass{article}

    \usepackage[final]{neurips_2022}

\usepackage[utf8]{inputenc} 
\usepackage[T1]{fontenc}    
\usepackage[unicode,colorlinks= true,allcolors=blue]{hyperref}       
\usepackage{url}            
\usepackage{booktabs,longtable}       
\usepackage{amsfonts}       
\usepackage{nicefrac}       
\usepackage{microtype}      
\usepackage{xcolor}         
\usepackage{enumitem}
\usepackage{lipsum}
\usepackage{amssymb}

\usepackage{subcaption}
\usepackage{placeins}
\usepackage{bm}
\usepackage{amsthm}
\usepackage{graphicx}
\usepackage{float}
\usepackage{array}
\usepackage{caption}
\usepackage[linesnumbered,algoruled,boxed,lined]{algorithm2e}

\usepackage{algpseudocode}

\algnewcommand{\Inputs}[1]{%
  \State \textbf{Inputs:}
  \Statex \hspace*{\algorithmicindent}\parbox[t]{.8\linewidth}{\raggedright #1}
}
\algnewcommand{\Initialize}[1]{%
  \State \textbf{Initialize:}
  \Statex \hspace*{\algorithmicindent}\parbox[t]{.8\linewidth}{\raggedright #1}
}
\SetAlFnt{\footnotesize}

\usepackage{mathtools}
 
\theoremstyle{definition}


\newcommand{\method}{\texttt{RLCG}}

\graphicspath{ {images/} }

\usepackage{wrapfig}

\title{A Deep Reinforcement Learning Framework for Column Generation}

\author{%
	Cheng Chi\\ University of Toronto\And
	Amine Mohamed Aboussalah\\ New York University \And
	Elias B. Khalil\thanks{Corresponding author: \href{mailto:khalil@mie.utoronto.ca}{\url{khalil@mie.utoronto.ca}}.} \\ University of Toronto\And
	Juyoung Wang\\ University of Toronto\And
	Zoha Sherkat-Masoumi \\ University of Toronto
}

\begin{document}

\maketitle

\begin{abstract}
Column Generation (CG) is an iterative algorithm for solving linear programs (LPs) with an extremely large number of variables (columns). CG is the workhorse for tackling large-scale \textit{integer} linear programs, which rely on CG to solve LP relaxations within a branch and price algorithm. Two canonical applications are the Cutting Stock Problem (CSP) and Vehicle Routing Problem with Time Windows (VRPTW). In VRPTW, for example, each binary variable represents the decision to include or exclude a \textit{route}, of which there are exponentially many; CG incrementally grows the subset of columns being used, ultimately converging to an optimal solution.  We propose~\method{}, the first Reinforcement Learning (RL) approach for CG. Unlike typical column selection rules which myopically select a column based on local information at each iteration, we treat CG as a sequential decision-making problem: the column selected in a given iteration affects subsequent column selections. This perspective lends itself to a Deep Reinforcement Learning approach that uses Graph Neural Networks (GNNs) to represent the variable-constraint structure in the LP of interest. 
We perform an extensive set of experiments using the publicly available BPPLIB benchmark for CSP and Solomon benchmark for VRPTW. 
\method{} converges faster and reduces the number of CG iterations by 22.4\% for CSP and 40.9\% for VRPTW on average compared to a commonly used greedy policy. Our code is available at \href{https://github.com/chichengmessi/reinforcement-learning-for-column-generation.git}{this link}.

 
\end{abstract}

\section{Introduction}
\label{sec:intro}

Machine Learning (ML) for Mathematical Optimization (MO) is a growing field with a wide range of recent studies that enhance optimization algorithms by embedding ML in them to replace human-designed heuristics~\citep{bengio2021machine}. Previous work has overwhelmingly focused on combinatorial or integer programming problems/algorithms for which all decision variables are explicitly given in advance. For example, in a knapsack problem, each binary variable represents the decision to include or exclude an item in the knapsack; even with hundreds of thousands of items, modern integer programming solvers or heuristics can assign values to some or all variables simultaneously without any memory issues. However, there are many optimization problems in which there are more decision variables than one could ever explicitly deal with. In VRPTW, for example, each binary variable represents the decision to include or exclude a \textit{route}, of which there are exponentially many. 

Column generation (CG) is an algorithm for solving linear programs (LPs) with a prohibitively large number of variables (columns)~\citep{DDS2006CG}. It leverages the insight that in an optimal solution to LP, only a few variables will be active. In each iteration of CG, a column with negative reduced cost is added to a Restricted Master Problem (RMP, where ``restricted" refers to the use of only a subset of columns) until no more columns with negative reduced cost are found. When that occurs, CG will have provably converged to an optimal solution to the LP. To solve an \textit{integer} linear program, CG can be used as an LP relaxation solver within the so-called \textit{branch and price} algorithm. In this work, we will focus on CG as an algorithm for solving large-scale linear programs; our conclusion section discusses the straightforward application of our method to the integer case.

A commonly used heuristic rule is to greedily select the column with the most negative reduced cost in each iteration. 
However, is this always the ``optimal" column to add if one is interested in converging in as few iterations as possible? Could ML-guided column selection, based on the structure of the current instance being solved, make better selection decisions that speed up the convergence of CG? These are the questions we tackle in this work. We adopt the standard view that one is interested in solving problem instances that have the same mathematical formulation but whose \textit{data} (objective function coefficients, constraint coefficients, num. of vehicles, etc.) differ and are drawn from the same distribution. In VRPTW, for example, the geographical locations of customers that must be served by the vehicles may vary, as do the corresponding service time windows. 




To tackle these questions, we propose a ML framework to accelerate the convergence of the CG algorithm. In particular, we utilize RL to select columns for LPs with many variables, where the state of the CG algorithm and the structure of the LP instance are encoded using GNNs that operate on a bipartite graph representation of the variable-constraint interactions, similar to~\citet{gasse2019exact}. Q-learning is used to derive a column selection policy that minimizes the total number of iterations through an appropriately designed reward function. Our premise is that by directly optimizing for convergence speed using RL on similar LP instances from the same problem (e.g., CSP, VRPTW), one can outperform traditional heuristics such as the greedy rule.
Our contributions can be summarized as follows:
\begin{enumerate}
    \item\noindent\textbf{CG as a sequential task:} We bring forth  the first integration of RL and CG which appropriately captures the sequential decision-making nature of CG algorithm. Our RL agent,~\method{}, learns a Q-function which takes the future rewards (namely, total number of iterations) into consideration and therefore can make better column selections at each step. In contrast with prior work on RL for mathematical optimization, ours is, to our knowledge, the first to tackle the very widely applicable class of problems with exponentially many variables.

    \item\noindent\textbf{Curricula for learning over a diverse set of instances:} In practice, instances of the same optimization problem may vary not only in their data, but also their size/complexity. For example, in the CSP, instances may vary in the number of rolls that must be cut. To enable efficient RL over a widely varying instance distribution, we show how a curriculum can be designed for a given dataset of instances with a minimal amount of domain knowledge.
    
    \item\noindent\textbf{Evidence of substantial improvements over existing myopic heuristics:} We evaluate~\method{} on the CSP, which is known as the representative problem in this domain~\citep{BAVC2004, DDS2006CG}, and the VRPTW, another widely studied and applied problem. We compare~\method{} with the  commonly used greedy column selection strategy and 
    an expensive, integer programming based one-step lookahead strategy described by \cite{morabit}. Our algorithm converges faster than both of these in terms of the number of iterations and total time. Our results show the value of considering CG as a sequential decision-making problem and optimizing the entire solving trajectory through RL.
\end{enumerate}

\section{Related Work}
\label{sec:related}
Work on the use of RL to guide iterative algorithms can be traced back to~\cite{zhang1995reinforcement}, who used RL for a scheduling problem. More recently,~\citet{khalil2017learning} and~\citet{bello2016neural} proposed deep RL for constructing solutions to graph optimization problems; the survey by~\citet{mazyavkina2021reinforcement} summarizes subsequent advances in this space. \citet{khalil2017learning} were the first to use GNNs in this setting, a line of work that has also grown substantially (e.g.,~\cite{cappart2021combinatorial}) including in integer programming~\citep{gasse2019exact}. Relatedly,~\citet{tang} apply RL to the cutting plane algorithm for integer linear programming, where a policy that selects ``good" Gomory cuts is derived using evolutionary strategies. Within the framework of Constraint Programming, \citet{cappart} present a deep RL approach for branching variable selection within an exact algorithm. Nonetheless, to our knowledge, the sequential decision-making perspective has not been leveraged in the context of column generation or integer/linear programming with many variables, a setting which is relevant to many practical applications such as resource allocation (e.g., CSP), routing problems (e.g., VRPTW), and airline crew scheduling~\citep{barnhart2003airline}. 

The closest work to ours is that of~\citet{morabit}. They formulate the column selection in each CG iteration as a column classification task where the label of each column (select or not) is given by an ``expert". This ``expert" performs a one-step lookahead to identify the column which maximally improves the LP value of the Restricted Master Problem (RMP) of the next iteration. This is done with an extremely time-consuming mixed-integer linear program (MILP).
The RMPs are encoded using bipartite graphs with columns nodes ($v$) and constraint nodes ($c$), where an edge between $v$ and $c$ in the graph indicates the contribution of a column $v$ to a constraint $c$. Each node of the graph is then annotated with additional useful information for column selection stored as node features. \citet{morabit} then use a GNN to imitate the node selection of the expert using supervision, similarly to what was done for branching by~\cite{gasse2019exact}.

In contrast, our work treats the column generation as a sequential decision-making problem and utilizes RL to select a column at each iteration of CG. Our GNN acts as a Q-function approximator that maximizes the total future expected reward. As such, our work focuses on directly reducing the total number of CG iterations, whereas~\citet{morabit} derive a classifier that  does not consider the interdependencies across iterations, treating them as independent. One approach to accelerate CG further is to add multiple columns per CG iteration \cite{Desaulniers}; we discuss this extension in the Conclusion.

\section{Preliminaries: Column Generation}
\label{sec31}

We will use the canonical Cutting Stock Problem (CSP) to describe the CG method as is typically done in textbooks on the topic~\citep{DDS2006CG}. CSP is a general resource allocation problem where the objective is to subdivide a given quantity of a resource into a number of predetermined allocations to meet certain demands so that the total usage of the resources (e.g., total number of size-fixed paper rolls) is minimized. Such minimization is achieved by finding a set of optimal divisions of each resource, or in other words, using a set of optimal cutting patterns to divide resources. Due to the CSP's combinatorial nature and its exponentially large set of possible patterns (variables), CG is used to solve the LP relaxation of CSP iteratively without explicitly enumerating all possible  patterns.


The CSP formulation and column generation algorithm we use is a common modification of \cite{Gilmore}. The set of all feasible patterns $\mathcal{P}$ that can be cut from a roll is defined as:
\begin{gather*}
\mathcal{P} =\bigg\{x_k \in \mathbb{N}^n  : \sum_{i =1}^{n} a_i x_{ik} \leq L, 
x_{ik} \geq 0 \quad \forall i \in \{1,2,\ldots,n\}, \forall k \in \{1,2,\ldots,|\mathcal{P}|\}   \bigg\}.
\end{gather*}
where each pattern $p \in \mathcal{P}$ is represented using a vector $x_k \in \mathbb{N}^n$. With $a_i$ being a possible cut length from a roll with length $L$, each element of $x_k$ specifies how many such cuts with length $a_i$ are included in pattern $p$. For instance, assume the length of the resource roll $L$ is 4m, so all possible $a_i$s are 1m, 2m, 3m and 4m, and then one possible cutting pattern $p$ is represented by  $x_k = (0,2,0,0)$. Let $\lambda_p$ be the number of times pattern $p$ is used. The formulation with $\lambda_p$ being a decision variable is:
\[
\min_{\lambda \in \mathbb{N}^{|\mathcal{P}|}} 
\bigg\{ \sum_{p \in \mathcal{P}} \lambda_p : \sum_{p \in \mathcal{P}} x_{ip}  \lambda_p = d_i \; \forall i \in \{1,2,\ldots,n\} \bigg\},
\]
where the objective function minimizes the total number of  patterns used, which is equivalent to minimizing the number of rolls used. The constraints ensure demand is met, while enforcing the integrality restriction on $\lambda_p$. 

This problem has an extremely large number of decision variables as $\mathcal{P}$ is exponentially large. Therefore, the problem is decomposed into the Restricted Master Problem (RMP) and the Sub-Problem (SP). The RMP is obtained by relaxing the integrality restrictions on $\lambda_p$ with an initial set $\tilde{\mathcal{P}}$ where $\tilde{\mathcal{P}}\subset \mathcal{P}$. The RMP formulation of the cutting stock problem is defined as follows: 
\begin{gather*}
\min_{\lambda \in \mathbb{N}^{|\tilde{\mathcal{P}}|}} 
\bigg\{ \sum_{p \in \tilde{\mathcal{P}}} \lambda_p : \sum_{p \in \tilde{\mathcal{P}}} x_{ip}  \lambda_p  = d_i \; \forall i \in \{1, \cdots, n\}, 
\lambda_p \geq 0 \; \forall  p \in \tilde{\mathcal{P}} \bigg\}.
\end{gather*}

    

    

The SP formulation of the provided cutting stock problem is defined as follows: 
\[
\min_{x \in \mathbb{N}^n} \left\{ \sum_{i = 1}^n \pi_i x_i : \sum_{i = 1}^n a_ix_i \leq L \right\},
\]
where $\pi_i$ is the dual value associated with the demand constraints. The SP is used to generate a pattern $p$ represented by vector $x \in \mathbb{N}^n$ with the most negative reduced cost, then adding that $p$ to $\tilde{\mathcal{P}}$ in the next iteration.
Below, we provide an  overview of the column generation algorithm, noting that our method will intervene in Step 4 to make a potentially non-greedy column selection decision:
\begin{enumerate}[nosep]
    \item Solve the RMP to obtain $\lambda^\star$ and $\bar{\pi}$;
    \item Update the SP objective function using $\bar{\pi}$;
    \item Solve SP to obtain $x_i^\star$;
    \item if $1-\sum_{i = 1}^{n} \bar{\pi_i} x_i^\star \leq 0 $, add the column associated with $x_i^\star$ and return to step 1, else Stop.
\end{enumerate}

\section{The~\method{} Framework}
\begin{figure}[htbp!]
\centering
      \includegraphics[width=\columnwidth]{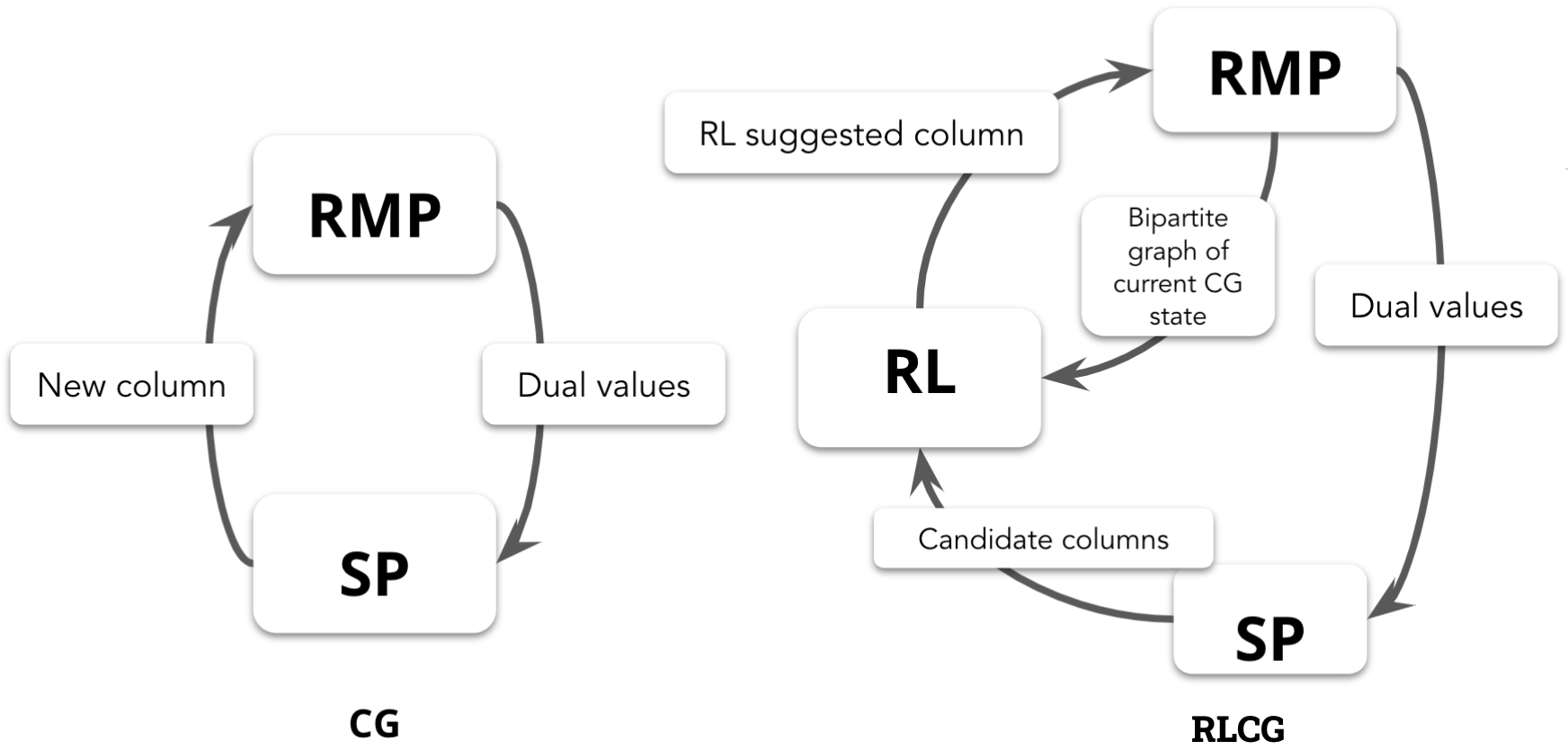}
      \caption{High-level comparison of standard Column Generation (CG) and~\method{}.}
      \label{fig:overview}
\end{figure} 


At a high level, our RL-aided column selection strategy,~\method{}, works as follows. We assume that the sub-problem (SP) is solved at each iteration and a set of near-optimal column candidates \( \mathcal{G} \) is returned, which is a general feature of optimization solvers such as Gurobi.
While the greedy CG algorithm, as described in Section \ref{sec31}, adds the single column with the most negative reduced cost from \( \mathcal{G} \) to the RMP of the next iteration, \method{} selects columns from \( \mathcal{G} \) according to the Q-function learned by the RL agent. The RL agent is fused within the CG loop and actively selects the column to be added to the next iteration using the information extracted from the current RMP and SP. An illustration comparing CG and~\method{} is provided in Figure \ref{fig:overview}.  

\subsection{Formulating CG as MDP}\label{sec42}
We formulate CG as a Markov decision process (MDP). As is customary, we use ($\mathcal{S}$, $\mathcal{A}$, $\mathcal{T}$, $r$, $\gamma$) to denote our MDP, where $\mathcal{S}$ is the state space, $\mathcal{A}$ the action space, $\mathcal{T}:\mathcal{S}\times\mathcal{S}\times\mathcal{A}\to[0,1]$, $(s',s,a)\mapsto\mathbb{P}(s'|s,a)$ the transition function, $r:\mathcal{S}\times\mathcal{A}\times\mathcal{A}\to\mathbb{R}$ the reward function, and $\gamma\in(0,1)$ the discount factor. We train the RL agent with experience replay \citep{mnih2013playing}. 

\subsubsection{State space $\mathcal{S}$}
\label{sec421}  

The state represents the information that the agent has about the environment. In~\method{}, the environment is the CG solution process corresponding to a given problem instance. As shown in Figure \ref{fig:overview}, the information passed to the 
RL agent is the bipartite graph of the current CG iteration from the RMP and the candidate columns from the SP. At each iteration, the RMP is an LP as shown in Section \ref{sec31}. As introduced in~\cite{gasse2019exact}, such an LP is encoded using a bipartite graph with two node types: column nodes $\mathcal{V}$ and constraint nodes $\mathcal{C}$. An edge $(v,c)$ exists between a node $v\in \mathcal{V}$ and a node $c\in \mathcal{C}$ if column $v$ contributes to constraint $c$. An example bipartite graph representation of the state is shown in left of Figure \ref{fig:state_trans} with column nodes shown on the left hand side, the constraint nodes on the right hand side, and the action nodes, which are the candidate columns returned from the SP, are shown in green (e.g., $v6,v7$).

To incorporate richer information about the current CG iteration, we include node features for both column nodes and constraint nodes (\textbf{vf} and  \textbf{cf} next to the nodes). We designed 9 column node features and 2 constraint node features in our environment based on our previous experience with CG and inspiration from \cite{morabit}. These node features are described in Appendix \ref{Features}.

As such, the state space $\mathcal{S}$ is the space of bipartite graphs representing all possible RMPs from the problem instances drawn from the distribution $\mathcal{D}$ with the node features given above. The bipartite graph shown in the left of Figure \ref{fig:state_trans} is a particular state $s$ in $\mathcal{S}$. As states are bipartite graphs with node features, it is natural to use a Graph Neural Network (GNN) as the Q-function approximator in our DQN agent. 
This bipartite graph representation not only encodes variable and constraint information in the RMP, but also interactions between variables and constraints through its edges.

\subsubsection{Actions, transition, and reward}\label{sec422}
\paragraph{Action space $\mathcal{A}$.}  

\begin{figure}[htb]
\centering
\includegraphics[scale=0.3]{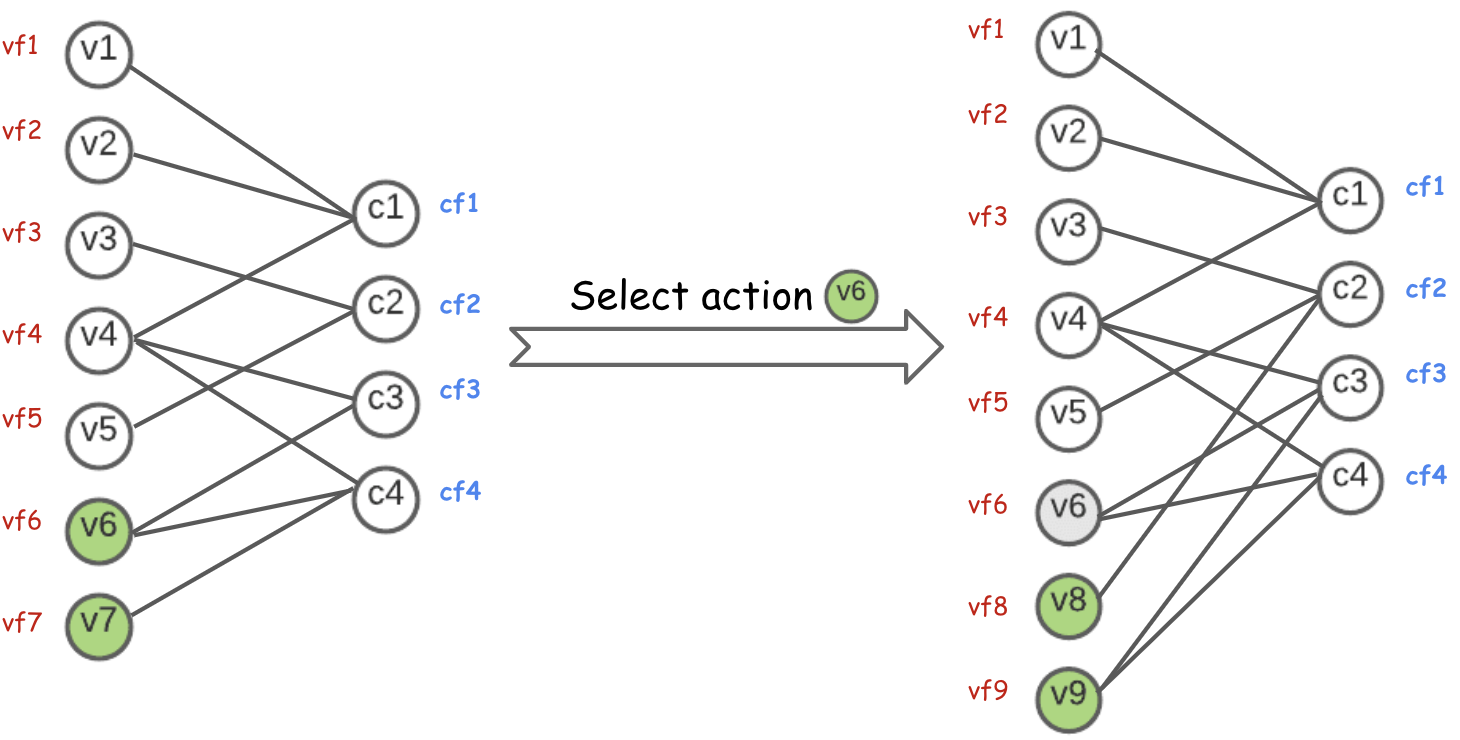}
\caption{State transition: Two green action nodes are considered, one of them is selected, transitioning to a new state.}
\label{fig:state_trans}
\end{figure}

As shown in Figure~\ref{fig:state_trans}, the RL agent selects one column to add to the RMP for the next iteration from the candidate set $\mathcal{G}$ returned from current iteration SP. Therefore, the action space $\mathcal{A}$ contains all possible candidate columns that can be generated from SPs; for example, the green nodes $v6$ and $v7$ in Figure \ref{fig:state_trans}. As the action space is discrete and the state space is continuous,  the GNN (Q-network) performs the operation of returning action values in the current state, i.e., $\hat{q}(s, a_1; w), \ldots, \hat{q}(s, a_m ;w)$.



\paragraph{Transition function $\mathcal{T}$.} 

Transitions are deterministic. After selecting an action from the current candidate set $\mathcal{G}$, the selected column enters the basis in the next RMP iteration. We then delete the action nodes that were not selected from the bipartite graph (current state), turn the selected action node into a column node, solve the new CG's RMP and SP, update all the features, and augment all the action nodes returned from the next SP iteration into the left-hand side of the graph, which results in a new bipartite graph state. Take the bipartite graph state shown in the left of Figure \ref{fig:state_trans} as an example, and assume action $v6$ is selected at this iteration. The transition occurs as  follows: $v6$ (in grey) becomes a column node, and candidate columns $v8$ and $v9$ returned by SP (in green) are added.

\paragraph{Reward function $\mathcal{R}$.} 
The reward function consists of two components: (1) the change in the RMP objective value, where a bigger decrease in value is preferred; (2) a unit penalty for each additional iteration. Together, they incentivize the RL agent to converge faster. 
The reward at time step $t$ is defined as:
\begin{gather}
\begin{gathered}
\label{eqn:reward}
r_t = \alpha\cdot\Bigg(\frac{\text{obj}_{t-1}-\text{obj}_{t}}{\text{obj}_{0}}\Bigg) 
- 1, 
\end{gathered}
\end{gather}
where obj$_0$ is the objective value of the RMP in the first CG iteration and is used to normalize $(\text{obj}_{t-1}-\text{obj}_{t})$ across instances of various sizes; $\alpha$ is a non-negative hyperparameter that weighs the normalized objective value change in the reward.

\subsection{~\method{} training and execution}
\label{sec43}

Algorithm \ref{alg:RL_CG} shows how a trained~\method{} agent is applied to solve a problem instance. The MDP components $(\mathcal{S},\mathcal{A},\mathcal{T},r,\gamma)$ used in this section are defined in Section \ref{sec42}.
Before starting the iterative optimization, the initialization steps 1--3 build the initial bipartite graph with computed node features and add an initial set of columns into basis (for instance, in CSP, we first add simple cutting patterns into the basis to initialize). Inside the while loop, a column is selected from the candidate set \( \mathcal{G} \) based on the best Q-value computed by the RL agent. Then, the RMP and the SP are updated in the same way as the traditional CG algorithm. The MDP model corresponds to extracting the MDP components from the current updated RMP and SP, which are discussed in Section \ref{sec421} and Section \ref{sec422}. Steps 7, 8 correspond to the deterministic state transition \( \mathcal{T} \) described in Section \ref{sec422}, and $S_{t+1}$ is the resulting state due to action $a^*_t$.

\begin{algorithm}[htb]
	\caption{~\method{} (RL-aided column generation algorithm)}
	\label{alg:RL_CG}
	
	\KwIn{Problem instance $p$ from distribution $\mathcal{D}$ \& trained Q-function.}
	\KwOut{Optimal solution}  
	$t = 0$;
	$\text{RMP}_0 = \text{Initialize}(p)$ \\
    Solve $\text{RMP}_0$ to get dual values;
    Use dual values to construct $\text{SP}_0$. \\
    $\langle S_0,A_0,T_0,R_0 \rangle = \text{MDP(RMP}_0\text{,SP}_0\text{)}$ \\
	\While{$\text{CG algorithm has not converged}$}{
	    $a^*_t = \text{arg}\max_{a_t \in A_t} \text{Q}(S_t,a_t) \quad \forall a_t \in A_t$  \;
	    Add variable $a^*_t$ to $\text{RMP}_t$ and get $\text{RMP}_{t+1}$ \\
	    Solve $\text{RMP}_{t+1}$ to get dual values;
	    Use dual values to build $\text{SP}_{t+1}$. \\
	    $\langle S_{t+1},A_{t+1},T_{t+1},R_{t+1} \rangle = \text{MDP(RMP}_{t+1}\text{,SP}_{t+1}\text{)}$ \\
	    $t = t+1$
	}
	Return optimal solutions from $\text{RMP}_\text{t}, \text{SP}_{t}$
\end{algorithm}

\textbf{Training.} The DQN algorithm with experience replay is used~\citep{mnih2013playing} with the typical mean squared loss between $Q(s_0,a)$ and $Q_{\text{target}}(s_0,a)$ $\forall$ $a$ (all green action nodes) in bipartite graph state $s_0$. $Q_{\text{target}}(s_0,a)$ is defined as $r+\gamma \cdot \max_{a} Q (s_1,a)$. A GNN is used as Q-function approximator. Training instances are sequenced based on a domain-specific curriculum described for each of CSP and VRPTW in Section~\ref{sec:Experiment}.

\section{Experimental Results}
\label{sec:Experiment}



\textbf{Baseline strategies and evaluation metrics.}
To assess the performance of~\method{}, we consider two baseline methods: the greedy column selection strategy (most negative reduced cost) and the MILP expert column selection of \cite{morabit} as described in Section~\ref{sec:related}. 
This expert provides an upper bound on the performance of the supervised learning approach of \cite{morabit}, as their ML model is approximately imitating the expert, i.e., it will never select better columns than the expert. We consider two standard evaluation metrics: (1) Number of iterations for CG to converge; (2) Time in seconds. For the latter, this includes GNN inference time and feature computations for~\method{}.
Our computing environment is described in Appendix Section~\ref{sec:computing}.

\subsection{Cutting Stock Problem (CSP)}
\label{sec:CSP}
\textbf{Dataset.}
We use BPPLIB~\cite{delorme2018bpplib}, a widely used collection of benchmark instances for binary packing and cutting stock problems, which includes a number of instances proven to be difficult to solve~\citep{delorme2020enhanced,wei2020new,martinovic2020improved}. BPPLIB contains instances of different sizes with the roll length $n$ varying from 50 to 1000 and the number of orders $m$ varying from 20 to 500. Our training set has instances with $n = {50, 100, 200}$. The remaining instances with $n = {50, 100, 200, 750}$ are split into a validation set and a testing set, with hard instances $n = 750$ only appearing in the testing set, as it is very expensive to solve such large instances during training. The detailed statistics of these three sets are listed in Table \ref{tab:datasets CSP} in Appendix \ref{Append: Train,Val,Test set numbers}. Note that our test set is more challenging than that used in training, allowing us to assess the agent's generalization ability. 

\textbf{Curriculum Design.}
In the RL context, curriculum learning \cite{narvekar2020curriculum} serves to sort the instances that a RL agent observes during the training process from easy to hard. In our experiments, we noticed that adopting the curriculum learning paradigm improves the learning of the RL agent and results in  better column selection and faster convergence. We train our~\method{} agent by feeding the instances in order of  increasing difficulty. For CSP, instances are sorted based on their roll length $n$ and number of orders $m$ to build a training curriculum (Table \ref{tab:learn} in Appendix C). The detailed comparison of training with and without the curriculum is deferred to Appendix~\ref{sec:curriculum}. In short, the former is crucial to convergence when the training set contains instances of varying difficulties.
The overall training process

\textbf{Hyperparameter tuning.}
We briefly describe the outcome of a large hyperparameter tuning effort which is described at length in Appendix~\ref{sec:hyperparam}. We tune the main hyperparameters, $\alpha$ in the reward function~\eqref{eqn:reward}, $\epsilon$ the exploration probability in DQN, $\gamma$ the discount factor, and the learning rate $lr$ in gradient descent. A grid with 81 possible configurations is explored by sampling 31 configurations, training the agent, and evaluating them on a validation set. The majority of the configurations resulted in agents that outperform the greedy strategy (see Appendix~\ref{sec:hyperparam} Figure~\ref{fig:hyper_val}). The best configuration is: $\alpha=300$,  $\epsilon = 0.05$, $\gamma = 0.9$ and $lr = 0.001$. For all the experiments conducted below, we use this configuration, both for CSP and VRPTW.

\begin{figure}[htb]
\begin{subfigure}{0.32\textwidth}
\includegraphics[width=\linewidth]{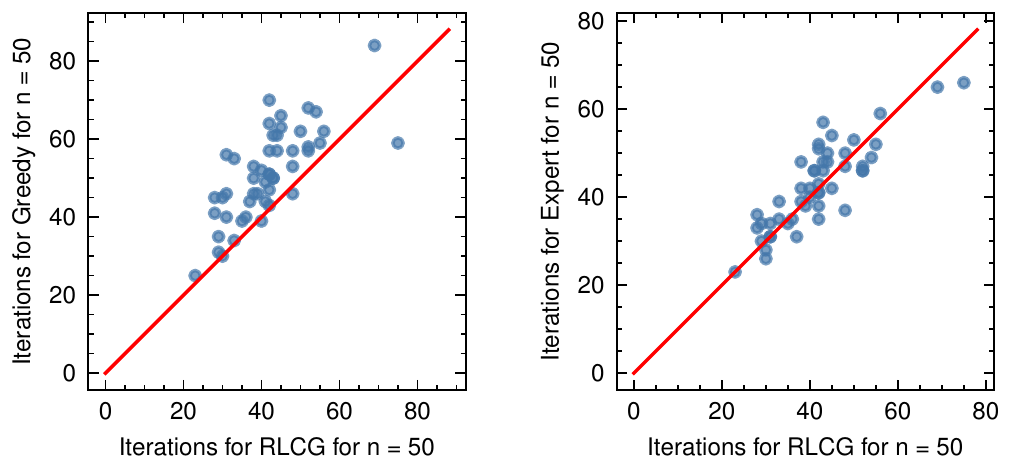}
\caption{CG iterations, $n=50$} \label{fig:7a}
\end{subfigure}\hspace*{\fill}
\begin{subfigure}{0.32\textwidth}
\includegraphics[width=\linewidth]{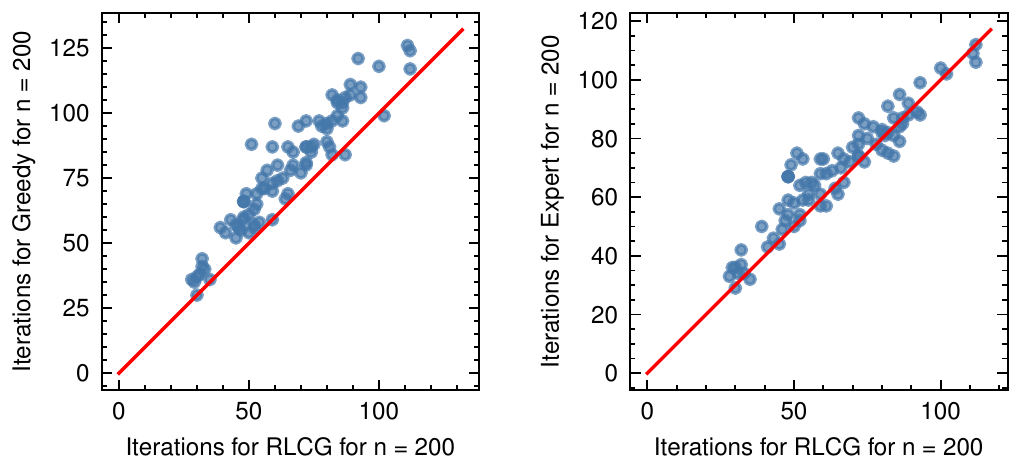}
\caption{CG iterations, $n=200$} \label{fig:7b}
\end{subfigure}\hspace*{\fill}
\begin{subfigure}{0.32\textwidth}
\includegraphics[width=\linewidth]{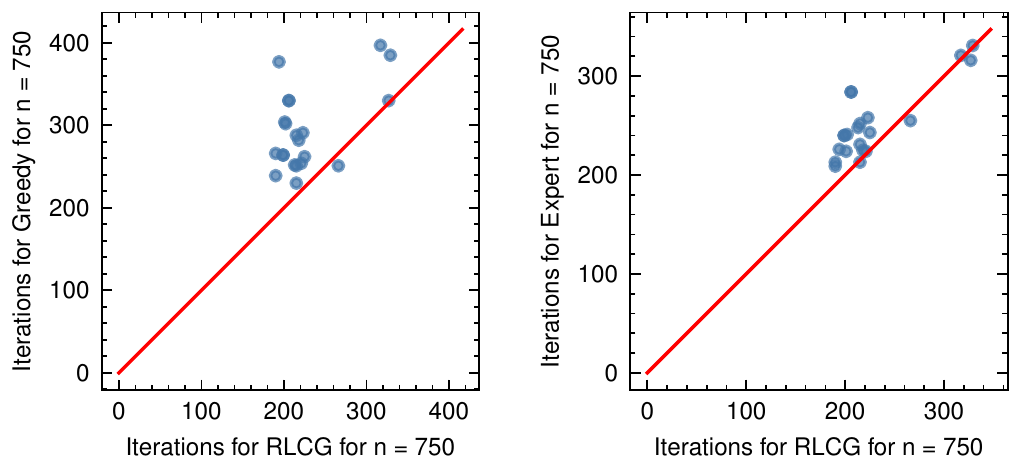}
\caption{CG iterations, $n=750$} \label{fig:7c}
\end{subfigure}

\medskip
\begin{subfigure}{0.32\textwidth}
\includegraphics[width=\linewidth]{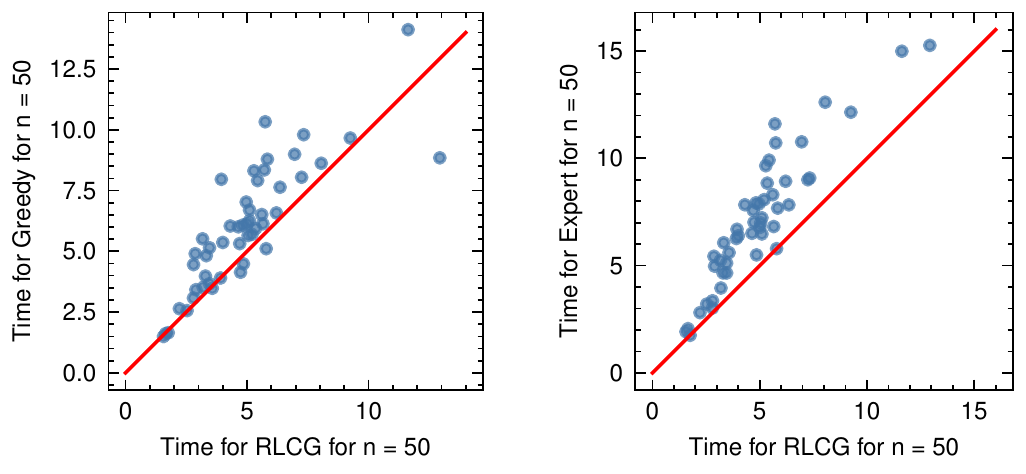}
\caption{Time in seconds, $n=50$} \label{fig:7d}
\end{subfigure}\hspace*{\fill}
\begin{subfigure}{0.32\textwidth}
\includegraphics[width=\linewidth]{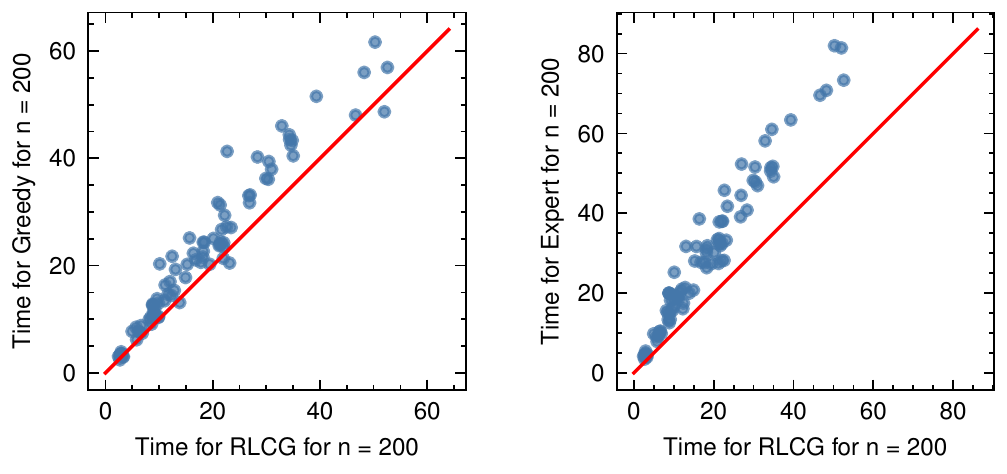}
\caption{Time in seconds, $n=200$} \label{fig:7e}
\end{subfigure}\hspace*{\fill}
\begin{subfigure}{0.32\textwidth}
\includegraphics[width=\linewidth]{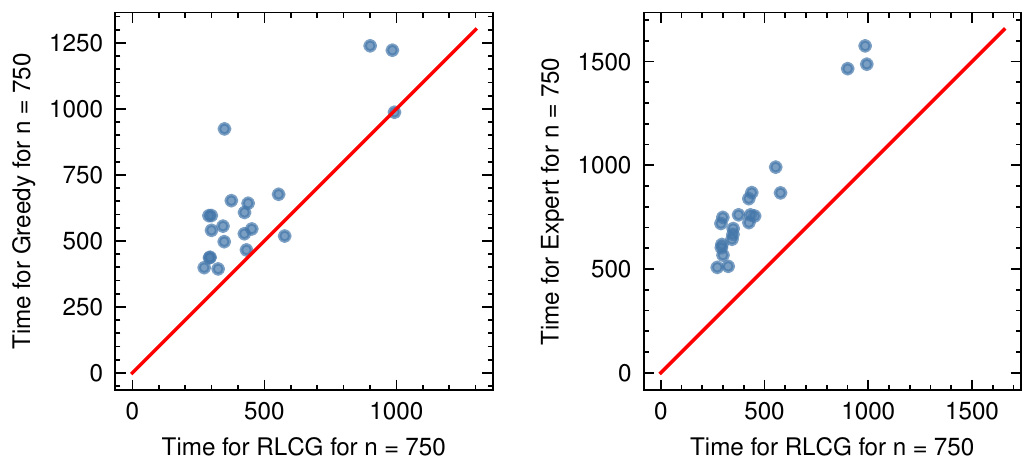}
\caption{Time in seconds, $n=750$} \label{fig:7f}
\end{subfigure}

\caption{\textbf{CSP:} Scatter plots of CG iterations (top row) and running time (bottom row) with~\method{} on the horizontal axis and Greedy (Expert) on the vertical axis, respectively, in each of three pairs of sub-figures. Each point represents a test instance of a given size $n$. Points above the diagonal indicate that~\method{} is faster than the competitor.} \label{fig:scatter}
\end{figure}

\textbf{Performance comparison.}
Figure \ref{fig:scatter} shows pair-wise comparisons between~\method{} and the two baselines in terms of the number of CG iterations (top row) and solving time (bottom row) on the testing instances. Each point in a plot corresponds to a test instance. Across all $n$ and for both baselines in Figure \ref{fig:scatter}, the majority of the points are above the diagonal line, which indicates that~\method{} outperforms the competing method w.r.t. the evaluation metric. Such a tendency becomes more pronounced as the CSP instances become larger (left to right).~\method{} performs better than the greedy column selection for all the roll lengths $n$ and the performance improvement becomes larger as $n$ grows. For CSP instances with $n=50$ and $n=200$, the MILP expert is able to maintain similar performance compared to RL. However, solving the MILP expert requires solving a MILP at each iteration with only one-step lookahead which is time consuming. For CSP instances with $n=750$,~\method{} begins to outperform the MILP expert due to its ability to take future rewards into consideration.

\begin{figure}[htb]
\begin{subfigure}{0.32\textwidth}
\includegraphics[width=\linewidth]{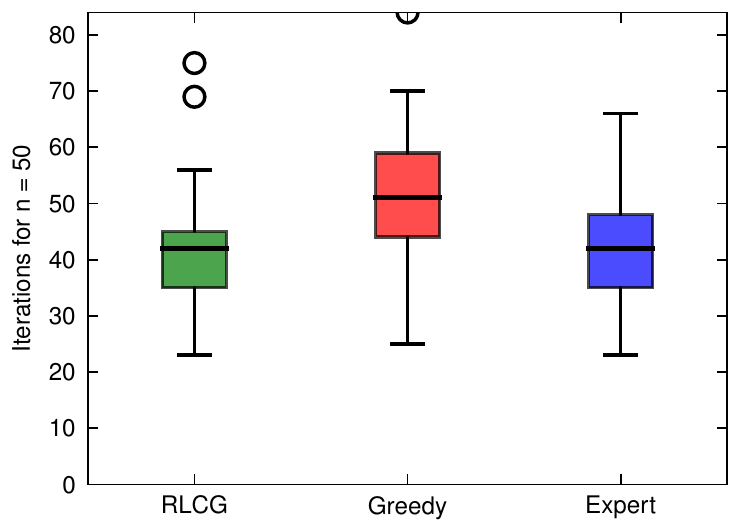}
\caption{CG iterations, $n=50$} \label{fig:a}
\end{subfigure}\hspace*{\fill}
\begin{subfigure}{0.32\textwidth}
\includegraphics[width=\linewidth]{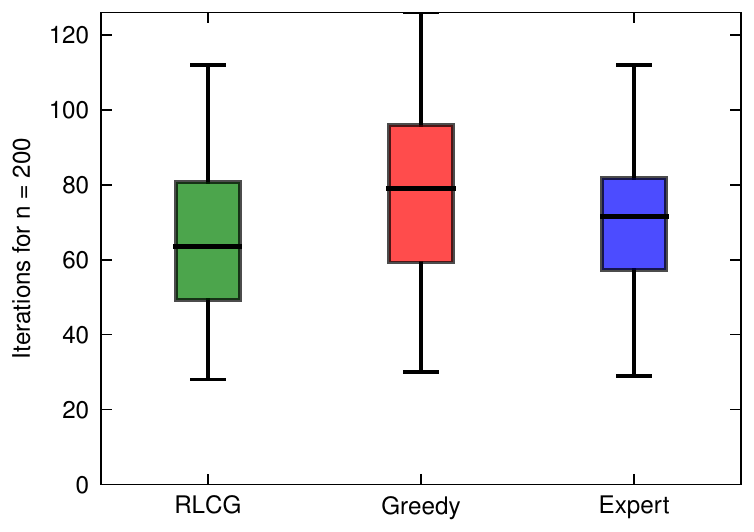}
\caption{CG iterations, $n=200$} \label{fig:a}
\end{subfigure}\hspace*{\fill}
\begin{subfigure}{0.32\textwidth}
\includegraphics[width=\linewidth]{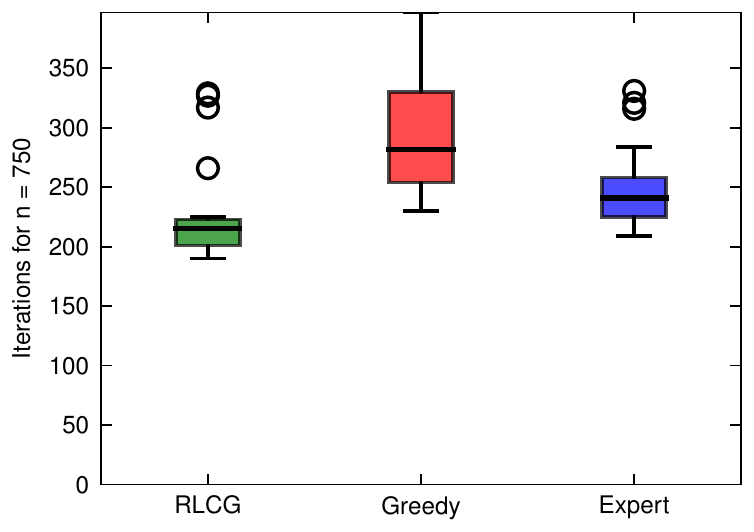}
\caption{CG iterations, $n=750$} \label{fig:b}
\end{subfigure}

\medskip
\begin{subfigure}{0.32\textwidth}
\includegraphics[width=\linewidth]{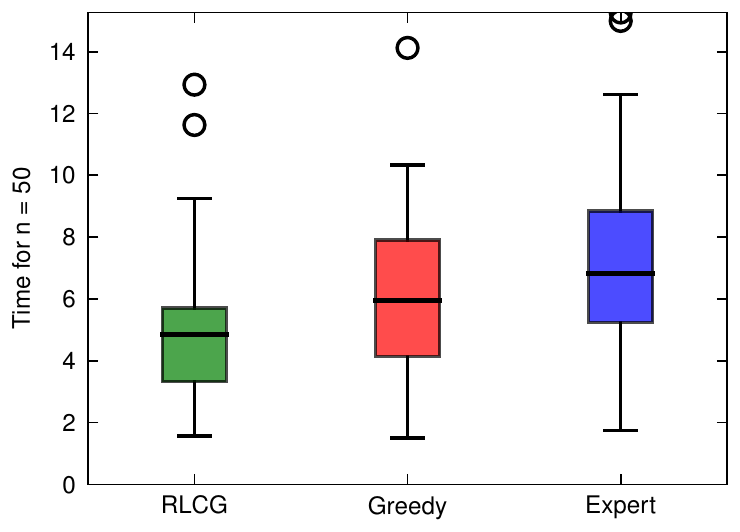}
\caption{Time in seconds, $n=50$} \label{fig:a}
\end{subfigure}\hspace*{\fill}
\begin{subfigure}{0.32\textwidth}
\includegraphics[width=\linewidth]{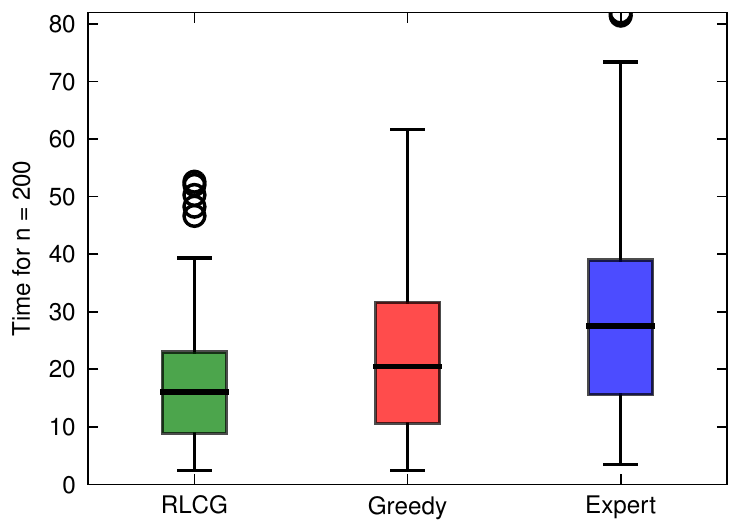}
\caption{Time in seconds, $n=200$} \label{fig:a}
\end{subfigure}\hspace*{\fill}
\begin{subfigure}{0.32\textwidth}
\includegraphics[width=\linewidth]{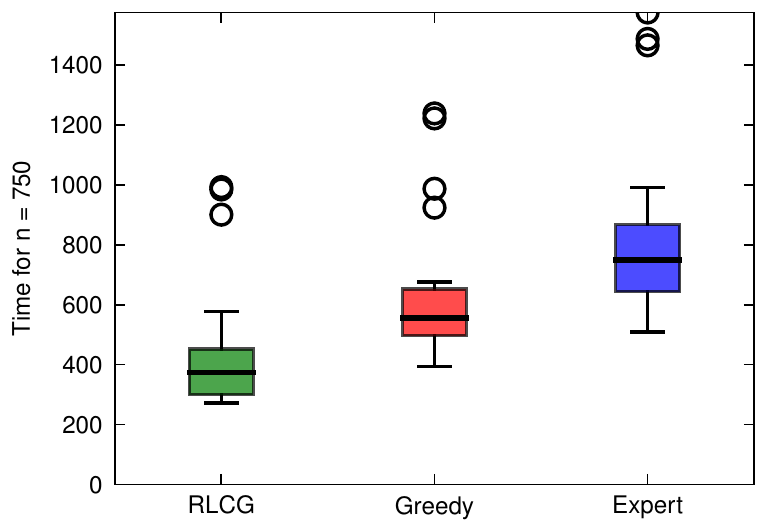}
\caption{Time in seconds, $n=750$} \label{fig:b}
\end{subfigure}
\caption{\textbf{CSP:} Box-plots of CG iterations (top row) and running time (bottom row) for the proposed method, RL (green), Greedy (red), and Expert (blue). Each box represents the distribution of the  iterations or running time for a given method and roll size $n$. Lower is better.} 
\label{fig:bar_plot}
\end{figure}

Notice that even though we did not train~\method{} using CSP instances with $n=750$, it was able to perform well on such instances, indicating strong generalization to harder problems. 
We also generate box plots to compare statistically the three methods shown in Figure \ref{fig:bar_plot}. Detailed statistics can be found in Table \ref{tab:statistics for CSP}. As shown by all the subplots in Figure \ref{fig:bar_plot}, especially for large instances in the subplots (c) and (f), the proposed~\method{} achieves statistically meaningful improvements on the CG convergence speed compared to both the greedy and the expert column selection methods. 

\begin{figure}[htb]
\begin{subfigure}{0.31\textwidth}
\includegraphics[width=\linewidth,scale=0.8]{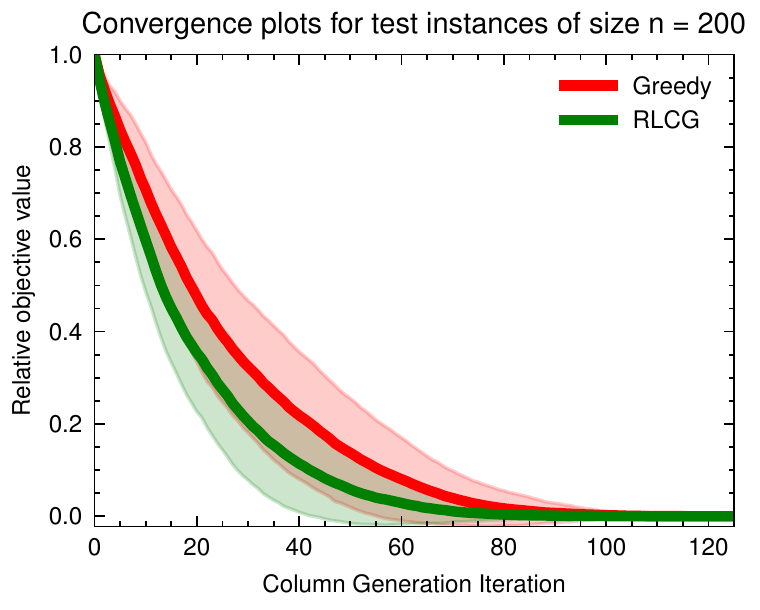}
\caption{\textbf{CSP}, $n=50$ test instances} \label{fig:a}
\end{subfigure}\hspace*{\fill}
\begin{subfigure}{0.31\textwidth}
\includegraphics[width=\linewidth,scale=0.8]{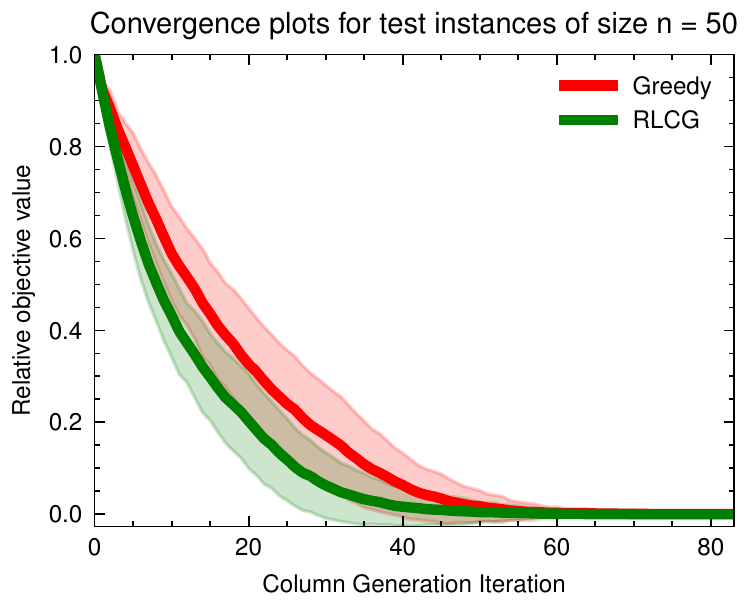}
\caption{\textbf{CSP}, $n=200$ test instances} \label{fig:a}
\end{subfigure}\hspace*{\fill}
\begin{subfigure}{0.31\textwidth}
\includegraphics[width=\linewidth,scale=0.8]{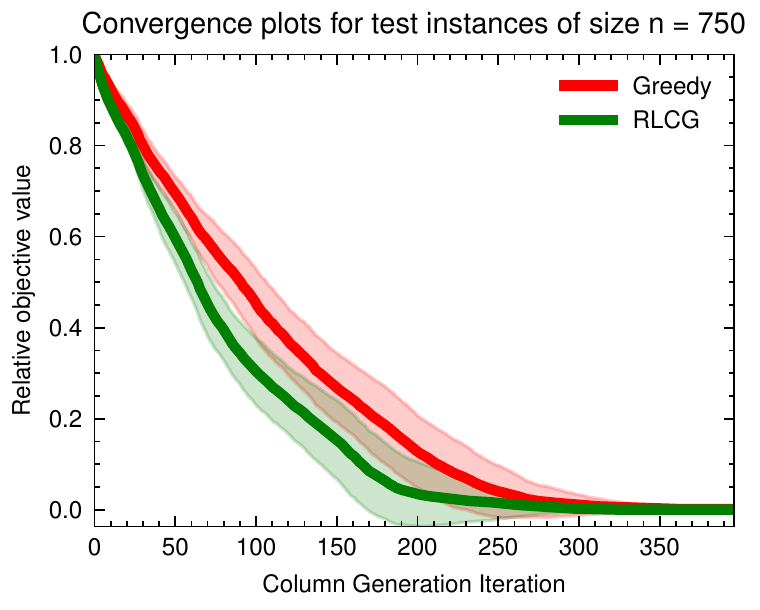}
\caption{\textbf{CSP}, $n=750$ test instances} \label{fig:a}
\end{subfigure}\hspace*{\fill}
\caption{CG convergence plots for \textbf{CSP}. The solid curves are the mean of the objective values for all instances and the shaded area shows $\pm 1$ standard deviation.}
\label{fig:convergence_CSP}
\end{figure}

\textbf{Convergence plots.}
In Figure~\ref{fig:convergence_CSP}, we visualize the CG solving trajectories for all test instances for $n=50,200,750$. We record the objective values of the RMP at each CG iteration for given method, then take the average over all instances. Note that we normalized the objective values to be in [0,1] before taking the average among instances. Since the MILP expert is by far the slowest, we restricted our attention to the two less expensive methods:~\method{} and Greedy. Looking at Figure~\ref{fig:convergence}, it is clear that not only does~\method{} terminate in fewer iterations and less time, but it also dominates Greedy throughout the CG iterations. In other words, if one had to terminate CG  earlier,~\method{} would result in a better (smaller) objective value as compared to what Greedy would achieve.

\subsection{Vehicle Routing Problem with Time Windows (VRPTW)}
\label{sec:VRPTW}

The VRPTW seeks a set of possible routes for a given number of vehicles to deliver goods to a group of geographically dispersed customers while minimizing the total travel costs. In the language of CG, each route is one column and there are exponentially many. Constraints of VRPTW include: vehicle capacity; a vehicle has to start from a depot and return to it; all customers should be served exactly once during their specified time windows. The detailed formulation of VRPTW is given by \cite{CordeauVRP2002}; implementation details are in the appendix.

\textbf{Dataset.}
We use the well-known Solomon benchmark \citep{solomon}. This dataset contains six different problem types (C1, C2, R1, R2, RC1, RC2), each of which has 8--12 instances with 50 customers. ``C" refers to customers which are geographically clustered, ``R" to randomly placed customers, ``RC" to a mixture. The ``1" and ``2" labels refer to narrow time windows/small vehicle capacity and large time windows/large vehicle capacity, respectively.
The difficulty levels of these sets are in order of C, R, RC. There are 56 instances in total in Solomon's dataset, and from each original Solomon instance, we can generate smaller instances by considering only the first $n<50$ customers.

We use instances from types C1, R1, RC1 for training. For the training set, we generated 80 smaller instances per type from the original Solomon's instances by only considering the first $n$ customers where $n$ is randomly sampled from 5--8, for a total of 240 training instances. For testing, two sets of instances are considered: 
(1) 60 small-size instances (number of customers $n$ within same range as training instances);
(2) 37 large-size instances (number of customers $n$ from 15--30). 
All the test instances are either generated from held-out Solomon instances (e.g., sets C2, R2, RC2) or generated from some training instance but with a much different $n$. 

\textbf{Curriculum Design.}
As in CSP, we designed a curriculum for VRPTW that sequences instances in order of difficulty: C1, R1, then RC1. This order is based on the fact that clustered instances have more structure that enables ``compact" routes for neighboring customers, whereas random instances require more complex routes. The mixed RC1 instances require reasoning about both types of customers simultaneously \citep{diff_VRP}.



\begin{figure}[htb]
\begin{subfigure}{0.5\columnwidth}
\includegraphics[width=.38\textwidth]{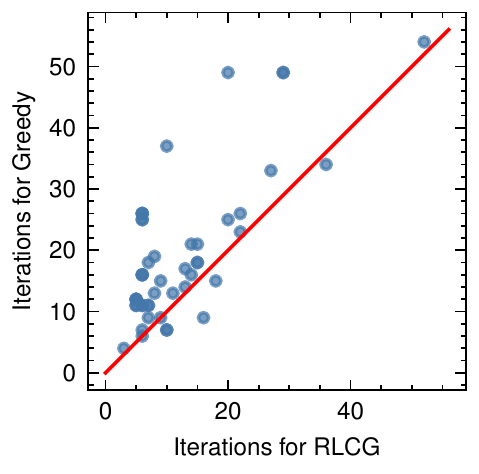}\quad
\includegraphics[width=.40\textwidth]{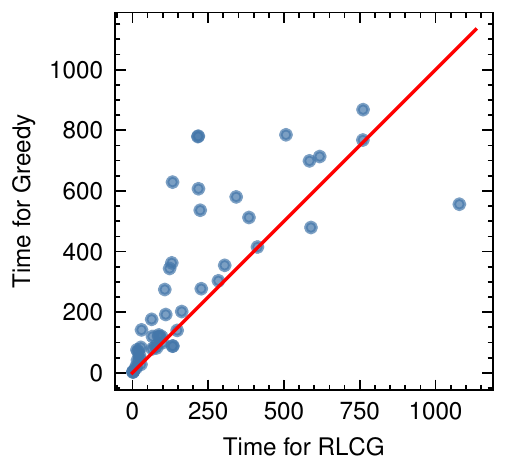} 
\caption{VRPTW with small instances} \label{fig:XXX}
\end{subfigure}
\hspace*{\fill}
\begin{subfigure}{0.5\columnwidth}
\includegraphics[width=.38\textwidth]{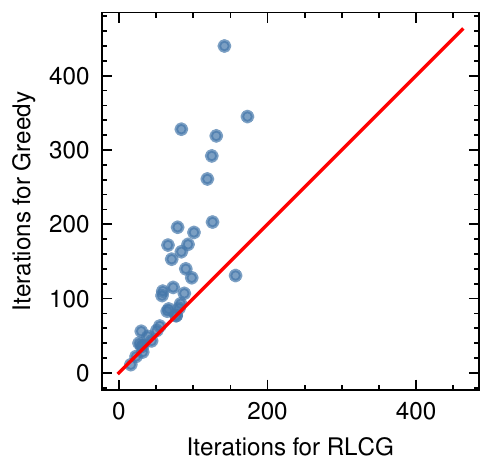} 
\quad
\includegraphics[width=.42\textwidth]{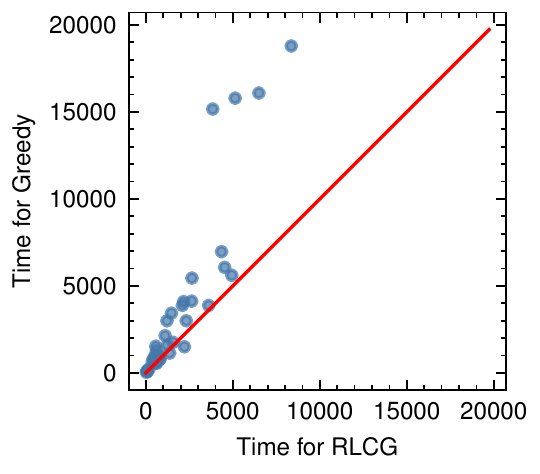} 
\caption{VRPTW with large instances} \label{fig:XXX}
\end{subfigure}
\caption{\textbf{VRPTW:} Scatter plots of CG iterations and running time with~\method{}. Each point represents a test instance of a given size $n$. Points above the diagonal indicate that~\method{} is faster than greedy.} 
\label{fig:scatterVRPTW}
\end{figure}

\textbf{Performance comparison.}
The hyperparemeters of the RL agent are chosen to be the same as the best set of hyperparameters found for CSP. Figure~\ref{fig:scatterVRPTW} shows similar trends to those seen for CSP:~\method{} converges in fewer iterations and less time than Greedy on most instances. 
\begin{figure}[htb]
\begin{subfigure}{0.45\textwidth}
\includegraphics[width=\linewidth,scale=0.8]{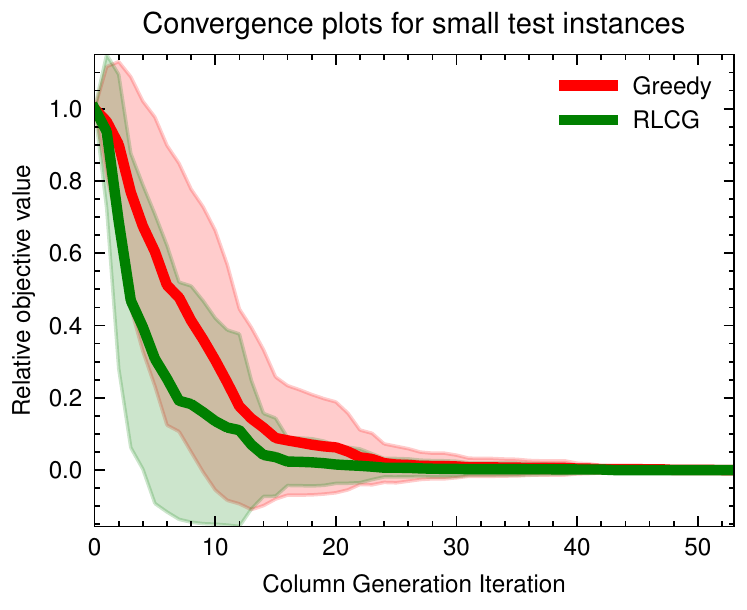}
\caption{\textbf{VRPTW}, Small test instances} \label{fig:a}
\end{subfigure}\hspace*{\fill}
\begin{subfigure}{0.45\textwidth}
\includegraphics[width=\linewidth,scale=0.8]{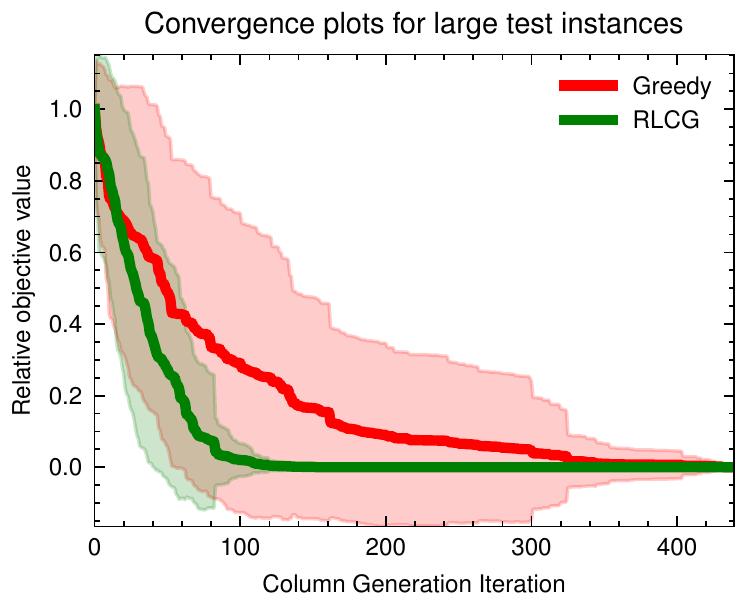}
\caption{\textbf{VRPTW}, Large test instances} \label{fig:a}
\end{subfigure}\hspace*{\fill}
\caption{CG convergence plots for \textbf{VRPTW}. The solid curves are the mean of the objective values for all instances and the shaded area shows $\pm 1$ standard deviation.}
\label{fig:convergence}
\end{figure}

This effect is more pronounced for the large VRPTW instances, a fact that can also be observed in the convergence plot of Figure~\ref{fig:convergence} (b):~\method{} converges in roughly 100 iterations compared to 300 iterations for Greedy, on average. Additional box plots are deferred to Appendix \ref{sec:VRPTW other plots}. Note that we do not compare to the Expert here given that it is much too slow and practically intractable.

\section{Conclusions \& Discussion of Limitations}
\label{sec:conclusion}

\method{} shows superior performance and better convergence both in terms of the number of iterations and time compared to the greedy column selection and the MILP expert strategy. In addition, our curriculum learning enables the agent to generalize well when facing harder test instances. Our experiments on two important large-scale LP families show that there is value in modeling CG as a sequential decision-making problem: taking the future impact of adding a column into account helps convergence.

However, our current work is restricted to adding only one column per CG iteration. Adding multiple columns per iteration can speed up convergence. However, this makes the RL action space exponential and thus finding the action with the largest Q-value becomes a combinatorial optimization problem. Recent work by~\citet{delarue2020reinforcement} addresses this problem, and thus our RL formulation can be expanded using the results of this paper to allow for multiple columns. Alternatively, policy gradient methods could be used instead of Q-learning. We believe this is an exciting next step but one that stretches beyond the limits of our paper, which is the first ever on RL for CG. A trained~\method{} agent can also be embedded within the branch and price algorithm for solving the integer-constrained versions of CSP/VRPTW and invoked to solve each LP relaxation. The speed-ups demonstrated herein would transfer to that setting, assuming that an appropriate dataset of training and validation instances can be collected.

\noindent{\textbf{Acknowledgments:}}
Khalil acknowledges support from the Scale AI Research Chair Program and an NSERC Discovery Grant.

\clearpage
\bibliographystyle{unsrtnat}
\bibliography{references}

\begin{thebibliography}{25}
\providecommand{\natexlab}[1]{#1}
\providecommand{\url}[1]{\texttt{#1}}
\expandafter\ifx\csname urlstyle\endcsname\relax
  \providecommand{\doi}[1]{doi: #1}\else
  \providecommand{\doi}{doi: \begingroup \urlstyle{rm}\Url}\fi

\bibitem[Bengio et~al.(2021)Bengio, Lodi, and Prouvost]{bengio2021machine}
Yoshua Bengio, Andrea Lodi, and Antoine Prouvost.
\newblock Machine learning for combinatorial optimization: a methodological
  tour d’horizon.
\newblock \emph{European Journal of Operational Research}, 290\penalty0
  (2):\penalty0 405--421, 2021.

\bibitem[Desaulniers et~al.(2006)Desaulniers, Desrosiers, and
  Marius]{DDS2006CG}
Guy Desaulniers, Jacques Desrosiers, and Solomon~M. Marius.
\newblock Column generation - a primer in column generation.
\newblock \emph{Springer}, page 1–32, 2006.

\bibitem[Gasse et~al.(2019)Gasse, Ch{\'e}telat, Ferroni, Charlin, and
  Lodi]{gasse2019exact}
Maxime Gasse, Didier Ch{\'e}telat, Nicola Ferroni, Laurent Charlin, and Andrea
  Lodi.
\newblock Exact combinatorial optimization with graph convolutional neural
  networks.
\newblock \emph{Advances in Neural Information Processing Systems}, 32, 2019.

\bibitem[Ben~Amor and Valerio~de Carvalho(2004)]{BAVC2004}
Hatem Ben~Amor and Jose Valerio~de Carvalho.
\newblock Cutting stock problems - in column generation.
\newblock \emph{Springer}, page 1–32, 2004.

\bibitem[Morabit et~al.(2021)Morabit, Desaulniers, and Lodi]{morabit}
Mouad Morabit, Guy Desaulniers, and Andrea Lodi.
\newblock Machine-learning--based column selection for column generation.
\newblock \emph{Transportation Science}, 55\penalty0 (4):\penalty0 815--831,
  2021.

\bibitem[Zhang and Dietterich(1995)]{zhang1995reinforcement}
Wei Zhang and Thomas~G Dietterich.
\newblock A reinforcement learning approach to job-shop scheduling.
\newblock In \emph{IJCAI}, volume~95, pages 1114--1120. Citeseer, 1995.

\bibitem[Dai et~al.(2017)Dai, Khalil, Zhang, Dilkina, and
  Song]{khalil2017learning}
Hanjun Dai, Elias Khalil, Yuyu Zhang, Bistra Dilkina, and Le~Song.
\newblock Learning combinatorial optimization algorithms over graphs.
\newblock \emph{Advances in neural information processing systems}, 30, 2017.

\bibitem[Bello et~al.(2016)Bello, Pham, Le, Norouzi, and
  Bengio]{bello2016neural}
Irwan Bello, Hieu Pham, Quoc~V Le, Mohammad Norouzi, and Samy Bengio.
\newblock Neural combinatorial optimization with reinforcement learning.
\newblock \emph{arXiv preprint arXiv:1611.09940}, 2016.

\bibitem[Mazyavkina et~al.(2021)Mazyavkina, Sviridov, Ivanov, and
  Burnaev]{mazyavkina2021reinforcement}
Nina Mazyavkina, Sergey Sviridov, Sergei Ivanov, and Evgeny Burnaev.
\newblock Reinforcement learning for combinatorial optimization: A survey.
\newblock \emph{Computers \& Operations Research}, 134:\penalty0 105400, 2021.

\bibitem[Cappart et~al.(2021)Cappart, Ch{\'e}telat, Khalil, Lodi, Morris, and
  Veli{\v{c}}kovi{\'c}]{cappart2021combinatorial}
Quentin Cappart, Didier Ch{\'e}telat, Elias Khalil, Andrea Lodi, Christopher
  Morris, and Petar Veli{\v{c}}kovi{\'c}.
\newblock Combinatorial optimization and reasoning with graph neural networks.
\newblock \emph{arXiv preprint arXiv:2102.09544}, 2021.

\bibitem[Tang et~al.(2020)Tang, Agrawal, and Faenza]{tang}
Yunhao Tang, Shipra Agrawal, and Yuri Faenza.
\newblock Reinforcement learning for integer programming: Learning to cut.
\newblock In \emph{International Conference on Machine Learning}, pages
  9367--9376. PMLR, 2020.

\bibitem[Cappart et~al.(2020)Cappart, Moisan, Rousseau, Pr{\'e}mont-Schwarz,
  and Cire]{cappart}
Quentin Cappart, Thierry Moisan, Louis-Martin Rousseau, Isabeau
  Pr{\'e}mont-Schwarz, and Andre Cire.
\newblock Combining reinforcement learning and constraint programming for
  combinatorial optimization.
\newblock \emph{arXiv preprint arXiv:2006.01610}, 2020.

\bibitem[Barnhart et~al.(2003)Barnhart, Cohn, Johnson, Klabjan, Nemhauser, and
  Vance]{barnhart2003airline}
Cynthia Barnhart, Amy~M Cohn, Ellis~L Johnson, Diego Klabjan, George~L
  Nemhauser, and Pamela~H Vance.
\newblock Airline crew scheduling.
\newblock In \emph{Handbook of transportation science}, pages 517--560.
  Springer, 2003.

\bibitem[Desaulniers et~al.(1999)Desaulniers, Desrosiers, and
  Solomon]{Desaulniers}
Guy Desaulniers, Jacques Desrosiers, and Marius Solomon.
\newblock Accelerating strategies in column generation methods for vehicle
  routing and crew scheduling problems.
\newblock 15, 01 1999.
\newblock \doi{10.1007/978-1-4615-1507-4_14}.

\bibitem[Gilmore and Gomory(1961)]{Gilmore}
P.~C. Gilmore and R.~E. Gomory.
\newblock A linear programming approach to the cutting-stock problem.
\newblock \emph{Operations research}, 9\penalty0 (6):\penalty0 849--859, 1961.
\newblock ISSN 0030-364X.
\newblock \doi{10.1287/opre.9.6.849}.
\newblock URL \url{https://doi.org/10.1287/opre.9.6.849}.

\bibitem[Mnih et~al.(2015)Mnih, Kavukcuoglu, Silver, Rusu, Veness, Bellemare,
  Graves, Riedmiller, Fidjeland, Ostrovski, et~al.]{mnih2013playing}
Volodymyr Mnih, Koray Kavukcuoglu, David Silver, Andrei~A Rusu, Joel Veness,
  Marc~G Bellemare, Alex Graves, Martin Riedmiller, Andreas~K Fidjeland, Georg
  Ostrovski, et~al.
\newblock Human-level control through deep reinforcement learning.
\newblock \emph{nature}, 518\penalty0 (7540):\penalty0 529--533, 2015.

\bibitem[Delorme et~al.(2018)Delorme, Iori, and Martello]{delorme2018bpplib}
Maxence Delorme, Manuel Iori, and Silvano Martello.
\newblock Bpplib: a library for bin packing and cutting stock problems.
\newblock \emph{Optimization Letters}, 12\penalty0 (2):\penalty0 235--250,
  2018.

\bibitem[Delorme and Iori(2020)]{delorme2020enhanced}
Maxence Delorme and Manuel Iori.
\newblock Enhanced pseudo-polynomial formulations for bin packing and cutting
  stock problems.
\newblock \emph{INFORMS Journal on Computing}, 32\penalty0 (1):\penalty0
  101--119, 2020.

\bibitem[Wei et~al.(2020)Wei, Luo, Baldacci, and Lim]{wei2020new}
Lijun Wei, Zhixing Luo, Roberto Baldacci, and Andrew Lim.
\newblock A new branch-and-price-and-cut algorithm for one-dimensional
  bin-packing problems.
\newblock \emph{INFORMS Journal on Computing}, 32\penalty0 (2):\penalty0
  428--443, 2020.

\bibitem[Martinovic et~al.(2020)Martinovic, Delorme, Iori, Scheithauer, and
  Strasdat]{martinovic2020improved}
John Martinovic, Maxence Delorme, Manuel Iori, Guntram Scheithauer, and Nico
  Strasdat.
\newblock Improved flow-based formulations for the skiving stock problem.
\newblock \emph{Computers \& Operations Research}, 113:\penalty0 104770, 2020.

\bibitem[Narvekar et~al.(2020)Narvekar, Peng, Leonetti, Sinapov, Taylor, and
  Stone]{narvekar2020curriculum}
Sanmit Narvekar, Bei Peng, Matteo Leonetti, Jivko Sinapov, Matthew~E Taylor,
  and Peter Stone.
\newblock Curriculum learning for reinforcement learning domains: A framework
  and survey.
\newblock \emph{arXiv preprint arXiv:2003.04960}, 2020.

\bibitem[Cordeau et~al.(2002)Cordeau, Desaulniers, Desrosiers, and
  Soumis]{CordeauVRP2002}
J.-F. Cordeau, G.~Desaulniers, J.~Desrosiers, and F.~Soumis.
\newblock The vrp with time windows.
\newblock \emph{SIAM Monographs on Discrete Mathematics and its Applications},
  pages 157--193, 2002.

\bibitem[Solomon(1987)]{solomon}
M.~Solomon.
\newblock Algorithms for the vehicle routing and scheduling problem with time
  window constraints.
\newblock \emph{Journal of Service Science and Management}, 35:\penalty0
  254--265, 1987.
\newblock URL
  \url{https://www.sintef.no/projectweb/top/vrptw/solomon-benchmark/}.

\bibitem[Desrochers(1992)]{diff_VRP}
Jacques Desrosiers; Marius Solomon;~Martin Desrochers.
\newblock A new optimization algorithm for the vehicle routing problem with
  time windows.
\newblock \emph{Operations Research}, 40\penalty0 (2):\penalty0 342--354, 1992.

\bibitem[Delarue et~al.(2020)Delarue, Anderson, and
  Tjandraatmadja]{delarue2020reinforcement}
Arthur Delarue, Ross Anderson, and Christian Tjandraatmadja.
\newblock Reinforcement learning with combinatorial actions: An application to
  vehicle routing.
\newblock \emph{Advances in Neural Information Processing Systems},
  33:\penalty0 609--620, 2020.

\end{thebibliography}
\clearpage
\section*{Checklist}

The checklist follows the references.  Please
read the checklist guidelines carefully for information on how to answer these
questions.  For each question, change the default \answerTODO{} to \answerYes{},
\answerNo{}, or \answerNA{}.  You are strongly encouraged to include a {\bf
justification to your answer}, either by referencing the appropriate section of
your paper or providing a brief inline description.  For example:
\begin{itemize}
  \item Did you include the license to the code and datasets? \answerYes{See \textbf{Dataset} and code link in appendix.}
  \item Did you include the license to the code and datasets? \answerNo{The code and the data are proprietary.}
  \item Did you include the license to the code and datasets? \answerNA{}
\end{itemize}
Please do not modify the questions and only use the provided macros for your
answers.  Note that the Checklist section does not count towards the page
limit.  In your paper, please delete this instructions block and only keep the
Checklist section heading above along with the questions/answers below.

\begin{enumerate}

\item For all authors...
\begin{enumerate}
  \item Do the main claims made in the abstract and introduction accurately reflect the paper's contributions and scope?
    \answerYes{} As discussed in Section \ref{sec:intro} and \ref{sec:related}, this work is the first application of RL on CG, and the benefits of using RLCG to enhance CG is demonstrated in Section \ref{sec:Experiment} and shown in Figure \ref{fig:scatter} \ref{fig:bar_plot}  \ref{fig:convergence} \ref{fig:scatterVRPTW}.Detailed improvement statistics can be verified in Table \ref{tab:statistics for VRP} \ref{tab:statistics for CSP}.
  \item Did you describe the limitations of your work?
    \answerYes{} Limitation is discussed at the end of Section \ref{sec:conclusion}. 
  \item Did you discuss any potential negative societal impacts of your work?
    \answerNA{}
  \item Have you read the ethics review guidelines and ensured that your paper conforms to them?
    \answerYes{}
\end{enumerate}

\item If you are including theoretical results...
\begin{enumerate}
  \item Did you state the full set of assumptions of all theoretical results?
    \answerNA{}
        \item Did you include complete proofs of all theoretical results?
    \answerNA{}
\end{enumerate}

\item If you ran experiments...
\begin{enumerate}
  \item Did you include the code, data, and instructions needed to reproduce the main experimental results (either in the supplemental material or as a URL)?
    \answerYes{}
    Code link is provided in Abstract. Data we use, training and testing procedure are discussed in detail in corresponding sections in \ref{sec:Experiment} and Appendix \ref{sec:curriculum} \ref{sec:computing} \ref{Append: Train,Val,Test set numbers}.  Hyperparameter settings for experiments in Appendix \ref{sec:hyperparam}.

    
  \item Did you specify all the training details (e.g., data splits, hyperparameters, how they were chosen)?
    \answerYes{}
        Datasets used, training testing split, training schedule generation are described set described in Section \ref{sec:CSP} for CSP and Section \ref{sec:VRPTW} for VRPTW. We conduct complete hyperparameter tuning for CSP in Appendix \ref{sec:hyperparam} and use the best hyperparameters setting for both CSP and VRPTW experiments.
        
        \item Did you report error bars (e.g., with respect to the random seed after running experiments multiple times)?
        
    \answerYes{} The solving process is deterministic so test on each instance once in the testing sets, however, we did normalize and average over all instance results and generate box plots in Figure \ref{fig:bar_plot} \ref{fig:box-VRP-small} and plot standard deviation bars in Figure  \ref{fig:convergence}.
        \item Did you include the total amount of compute and the type of resources used (e.g., type of GPUs, internal cluster, or cloud provider)? 
    \answerYes{}
        Computing environment in Appendix \ref{sec:computing};
    
\end{enumerate}

\item If you are using existing assets (e.g., code, data, models) or curating/releasing new assets...
\begin{enumerate}
  \item If your work uses existing assets, did you cite the creators?
    \answerNA{}
  \item Did you mention the license of the assets?
    \answerNA{}
  \item Did you include any new assets either in the supplemental material or as a URL?
    \answerNA{}
  \item Did you discuss whether and how consent was obtained from people whose data you're using/curating?
    \answerNA{}
  \item Did you discuss whether the data you are using/curating contains personally identifiable information or offensive content?
    \answerNA{}
\end{enumerate}

\item If you used crowdsourcing or conducted research with human subjects...
\begin{enumerate}
  \item Did you include the full text of instructions given to participants and screenshots, if applicable?
    \answerNA{}
  \item Did you describe any potential participant risks, with links to Institutional Review Board (IRB) approvals, if applicable?
    \answerNA{}
  \item Did you include the estimated hourly wage paid to participants and the total amount spent on participant compensation?
    \answerNA{}
\end{enumerate}

\end{enumerate}

\newpage
\appendix


\section{Curriculum learning}
\label{sec:curriculum}

\begin{figure*}[h!] 
    \centering
    \subfloat[curriculum training]{%
        \includegraphics[width=0.4\linewidth]{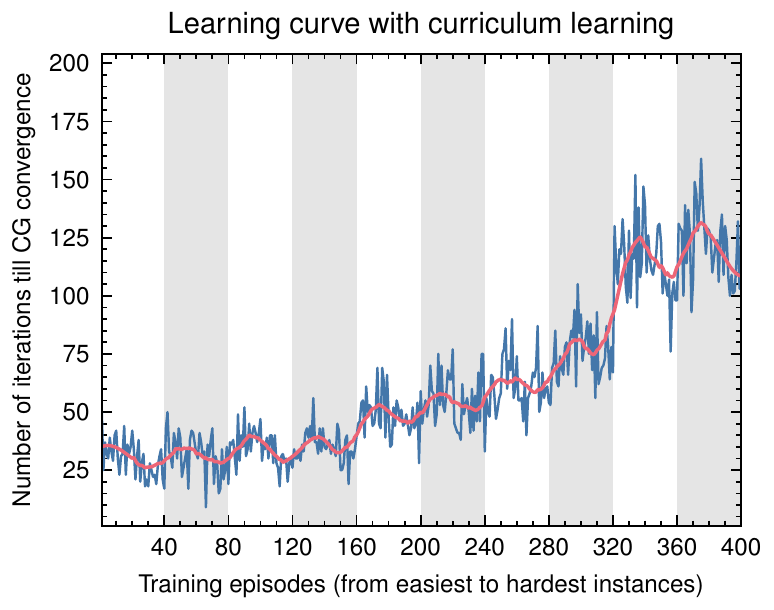}%
        \label{fig:a}%
        }%
  \quad
    \subfloat[No curriculum training]{%
        \includegraphics[width=0.4\linewidth]{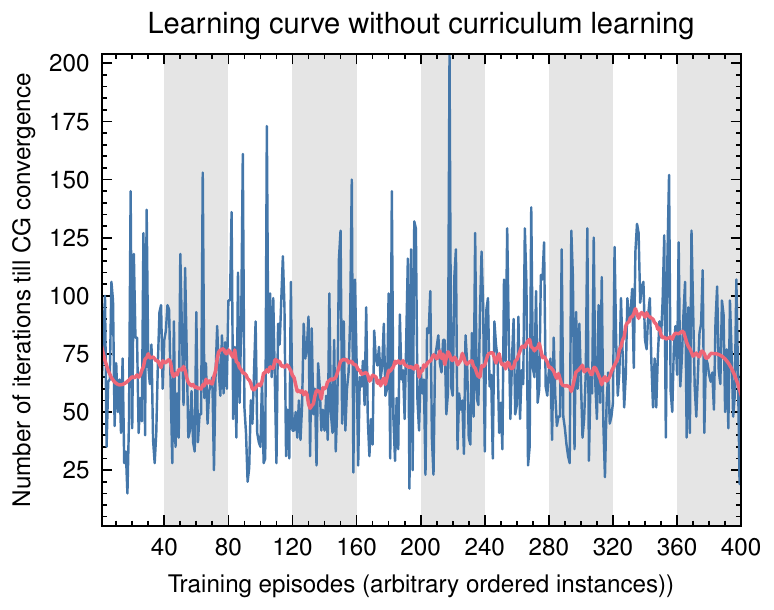}%
        \label{fig:b}%
        }%
    \caption{(a) shows the training process with such curriculum and Figure (b) shows the training process without a curriculum. In both plots, the x axis is the training episode, and each episode is solving one instance (instances are ordered for (a) and randomized for (b)) while the y axis shows the total number of RL guided CG iterations until that particular instance is solved.} \label{fig:no_s_vs_s}
\end{figure*}

Figure \ref{fig:no_s_vs_s} above compares the training trajectories between training with 400 CSP instances following the sequence provided in Table \ref{tab:learn} and training with the same 400 CSP instances randomly ordered. In (a), as every 40 instances we increase the CSP difficulty. Although there is an upward trend in the training curve, however, within each instance difficulty setting (fixed n and m), there is a downward trend showing a sign of learning. In contrast, there is no clear sign of learning in (b). Therefore, for all the experiments shown in this paper, the~\method{} model is trained using a curriculum. 

\begin{figure*}[h!] 
    \centering{%
        \includegraphics[width=0.4\linewidth]{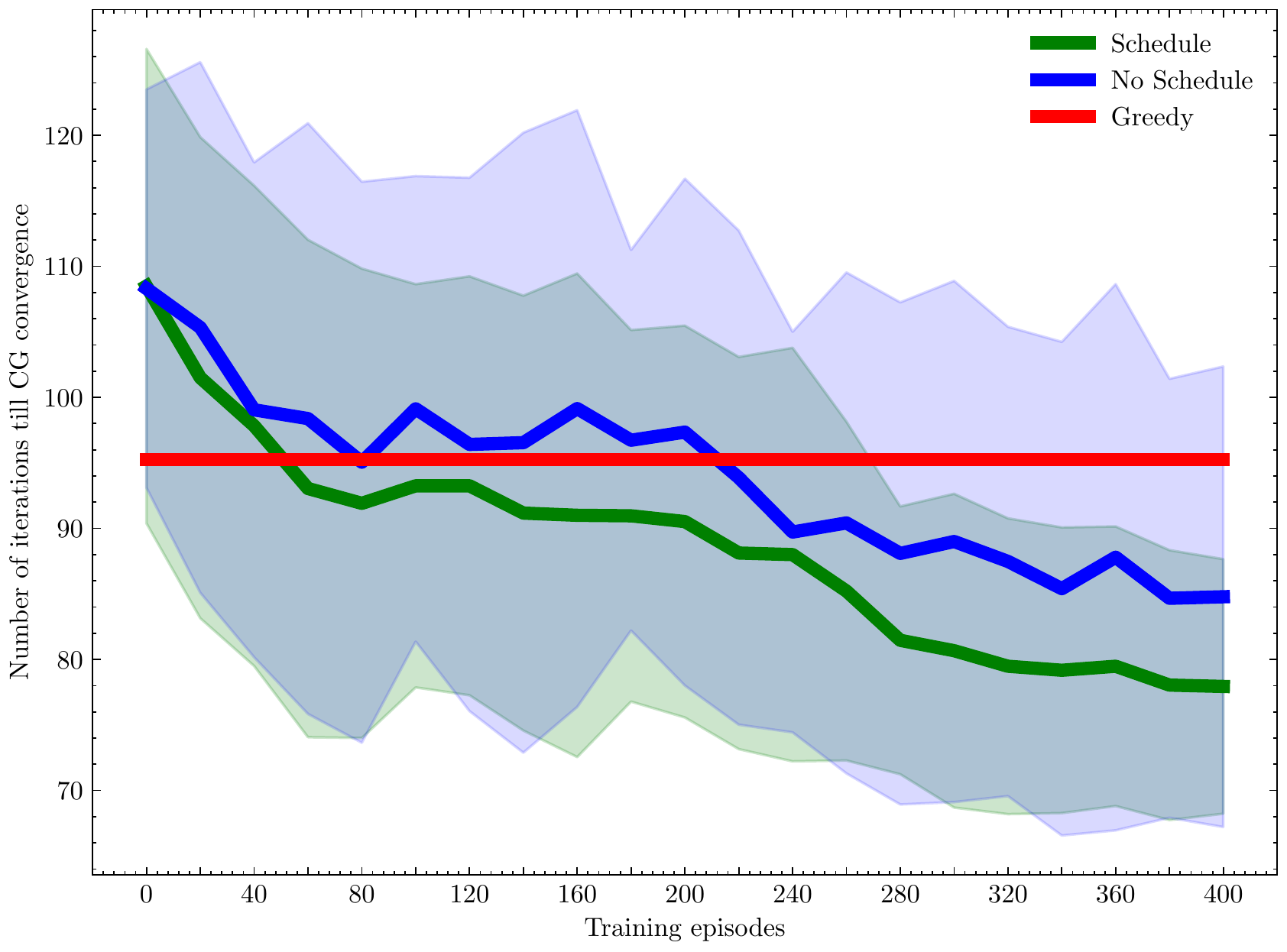}}
    \caption{Green curve shows the mean of convergence iterations for the current RLCG model with schedule training solving all the validation instances, and plus minus one standard deviation is shown in shaded area. Blue curve shows validation result of model trained without schedule. The red horizontal line shows average convergence iterations for validation instances using greedy strategy} \label{fig:val_train}
\end{figure*}

We also visualize the training process of RLCG for CSP using a validation set with 30 instances. The validation set detail is in Appendix \ref{Append: Train,Val,Test set numbers}. For every 20 training episodes, we stop the training process and validate the current models (with schedule training and without schedule training) with the validation set. The result is shown in Figure \ref{fig:val_train}.

\section{Hyperparameter tuning}
\label{sec:hyperparam}
Table \ref{tab:RL param} shows RL related parameters.

{\renewcommand{\arraystretch}{1.2}
\begin{table}[htb]\small
\caption{RL Parameters}
\label{tab:RL param}
\begin{center}
\begin{tabular}{|>{\raggedright}p{20mm}|p{55mm}|}
\hline
Parameter &  Value\\
\hline
state normalization & features normalization is \textit{MinMaxScalar} from sklearn \\
\hline
step penalty & For each iteration~\method{} takes before the current CSP instance is solved, penalize each step by 1 in the reward design \\
\hline
reward & there are two settings $\alpha=5$ and $\beta=1$, $\alpha=5$ and $\beta=0$, step penalty is always 1 \\
\hline
action & solution pool 10 from Gurobi solver for SP, which means at each iteration our action space contains 10 columns.\\
\hline
\end{tabular}
\end{center}
\end{table} 
}

In Table \ref{tab:GNN param}, we provide GNN parameters that we used.

\begin{figure}[htb]
    \centering
    \includegraphics[scale=1]{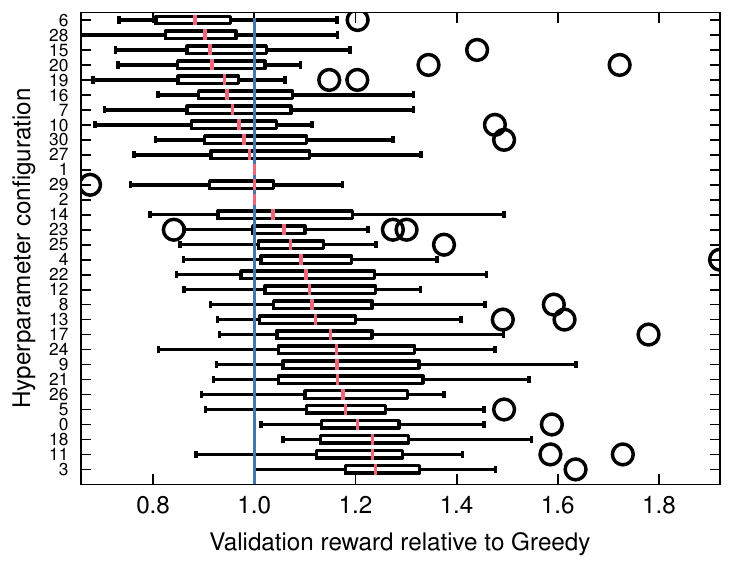}
    \caption{Hyperparameter sensitivity: the vertical axis list all the models we evaluated by its index. The detailed hyperparameter configurations each index refers to are listed in Appendix \ref{sec:hyperparam} Table \ref{fig:hyper_val}. The horizontal axis shows the ratio of~\method{} model solving iterations to greedy column selection solving iterations. The blue vertical line shows the ratio threshold indicating same performance as greedy. The best to worst models are ordered from bottom to top.}
    \label{fig:hyper_val}
\end{figure}

{\renewcommand{\arraystretch}{1.2}
\begin{table}[htb]\small
\caption{GNN Parameters}
\label{tab:GNN param}
\begin{center}
\begin{tabular}{|p{3cm}|p{4cm}| }
\hline
Parameter &  Value\\
\hline
Optimizer & Adam\\
\hline
Network Structure & refer to code for details\\
\hline
Batch Size & 32\\
\hline
\end{tabular}
\end{center}
\end{table}
}
We conduct hyperparameters tuning and sensitivity analysis using a validation set over: the parameter $\alpha$ used in equation \ref{eqn:reward} to weight the change in the normalized objective value used in our reward function, the exploration parameter of the RL agent $\epsilon$, the discount factor $\gamma$, and the learning rate $lr$. All other hyperparameters and their values are listed in Table \ref{tab:RL param} and Table \ref{tab:GNN param}.  
The values we consider for each hyperparameter are the following: $\alpha \in$ \{0, 100, 300\}, $\epsilon \in$ \{0.01, 0.05, 0.2\}, $\gamma \in$ \{0.9, 0.95, 0.99\} and $lr$ $\in$ \{0.01, 1e-3, 3e-4\}. We choose the value for $\alpha$ $\in$ \{0, 100, 300\} because when $\alpha = 0$, we place no weight on decreasing the objective value, otherwise, $\alpha=100, 300$ will bring the normalized change in the objective values into similar scale as the step penalty. Therefore, the search space for hyperparameters is defined as the Cartesian product between all these sets of different hyperparameters possible values, which gives us 81 configurations, and we randomly select 31 configurations out of them. Then we train 31~\method{} models corresponding to selected 31 configurations, and we evaluate these models using our validation set. The validation metric or reward is defined as the ratio of the total number of iterations~\method{} takes to solve each CSP instance divided by the total number of iterations greedy takes. For each model, we compute such ratio for all the instances in validation set, and we generate the following box plot shown in Figure \ref{fig:hyper_val} showing the validation metric for all the models. Among the 31 configurations we tested, the majority of the models were able to outperform greedy strategy on the instances in the validation set. The best configuration is model 3:  $\alpha=300$,  $\epsilon = 0.05$, $\gamma = 0.9$ and $lr = 0.001$. For all the results reported in this paper, we use this configuration. This includes VRPTW, although no direct hyperparameter tuning has been done for this problem class. To assess the sensitivity of the RL training with respect to randomness such as the GNN initialization and the exploration in RL, we compare the average validation reward relative to greedy for the selected model across five random seeds. The average varies between 1.23 and 1.25, indicating little to no sensitivity.

In Table \ref{tab:hyper_index}, we provide detailed configurations of each model index for our analysis of hyperparameters as well as their detailed validation results.

\newpage
\begin{table}[htb]\small
\begin{tabular}{llrrr}
\toprule
Model index &     Hyperparameters config &  iterations mean &  iterations median &  iterations std \\
\midrule
    Model 0 &    (100, 0.2, 0.9, 0.0003) &             1.22 &               1.20 &            0.14 \\
    Model 1 &    (100, 0.05, 0.99, 0.01) &             1.00 &               1.00 &            0.00 \\
    Model 2 &      (100, 0.2, 0.9, 0.01) &             1.00 &               1.00 &            0.00 \\
    Model 3 &    (300, 0.05, 0.9, 0.001) &             1.25 &               1.24 &            0.13 \\
    Model 4 &   (300, 0.05, 0.99, 0.001) &             1.12 &               1.09 &            0.19 \\
    Model 5 &   (300, 0.05, 0.95, 0.001) &             1.18 &               1.18 &            0.14 \\
    Model 6 &     (0, 0.05, 0.9, 0.0003) &             0.91 &               0.88 &            0.13 \\
    Model 7 &     (0, 0.2, 0.95, 0.0003) &             0.98 &               0.96 &            0.14 \\
    Model 8 &   (100, 0.05, 0.9, 0.0003) &             1.14 &               1.11 &            0.16 \\
    Model 9 &     (100, 0.2, 0.9, 0.001) &             1.19 &               1.16 &            0.18 \\
   Model 10 &    (0, 0.05, 0.95, 0.0003) &             0.97 &               0.97 &            0.14 \\
   Model 11 &   (300, 0.01, 0.95, 0.001) &             1.22 &               1.23 &            0.16 \\
   Model 12 &    (300, 0.05, 0.9, 0.001) &             1.11 &               1.11 &            0.12 \\
   Model 13 &  (100, 0.05, 0.95, 0.0003) &             1.15 &               1.12 &            0.16 \\
   Model 14 &     (100, 0.2, 0.9, 0.001) &             1.07 &               1.04 &            0.17 \\
   Model 15 &    (100, 0.05, 0.9, 0.001) &             0.95 &               0.91 &            0.15 \\
   Model 16 &      (0, 0.2, 0.99, 0.001) &             0.99 &               0.95 &            0.13 \\
   Model 17 &    (300, 0.2, 0.95, 0.001) &             1.18 &               1.15 &            0.17 \\
   Model 18 &    (300, 0.2, 0.9, 0.0003) &             1.24 &               1.23 &            0.13 \\
   Model 19 &     (0, 0.05, 0.9, 0.0003) &             0.92 &               0.94 &            0.12 \\
   Model 20 &       (0, 0.2, 0.9, 0.001) &             0.95 &               0.92 &            0.19 \\
   Model 21 &  (300, 0.05, 0.99, 0.0003) &             1.18 &               1.16 &            0.16 \\
   Model 22 &    (300, 0.2, 0.99, 0.001) &             1.11 &               1.10 &            0.16 \\
   Model 23 &      (100, 0.2, 0.9, 0.01) &             1.06 &               1.06 &            0.10 \\
   Model 24 &    (300, 0.05, 0.9, 0.001) &             1.16 &               1.16 &            0.17 \\
   Model 25 &    (100, 0.2, 0.95, 0.001) &             1.07 &               1.07 &            0.12 \\
   Model 26 &    (100, 0.01, 0.9, 0.001) &             1.19 &               1.18 &            0.13 \\
   Model 27 &       (0, 0.2, 0.9, 0.001) &             1.00 &               0.99 &            0.13 \\
   Model 28 &       (0, 0.2, 0.9, 0.001) &             0.90 &               0.90 &            0.11 \\
   Model 29 &     (0, 0.05, 0.95, 0.001) &             0.97 &               1.00 &            0.11 \\
   Model 30 &  (100, 0.05, 0.95, 0.0003) &             1.01 &               0.98 &            0.15 \\
\bottomrule
\end{tabular}
\caption{Evaluated models' configurations and validation performances}
\label{tab:hyper_index}
\end{table}

\newpage
\section{Computing environment}
\label{sec:computing}
To implement GNN, we use Tensorflow 2.7.0. To solve both RMP and SP optimization problems in CG, we use Gurobi 9.5.0. For the training using 400 instances schedule for CSP and 240 instances schedule for VRPTW, the training takes around 8-10 hours CPU time using the following CPU settings: Intel 2.30 Ghz, 2 CPU Cores, Intel (R) Xeon(R), Haswell CPU family. 

\section{Graph neural network for bipartite graph}
\label{sec:GNN_for_bipartite}

Graph Neural Network (GNN) has been successfully applied to many different machine learning tasks with graph structured data. GNN includes a message passing method where the features of each node for each node pass to other neighbouring nodes through learned transformations to generate aggregated information, and such information can be used for node classification, edge selection so on so forth. Due to the effectiveness of GNN to utilize graphical structure of the data, catch node-level dependencies, and has permutation invariant properties, GNN is an appropriate method for performing node(column) selection task in the present column generation problem. 

For this study we are encoding the information for each iteration of CG using a bipartite graph $G=(E,V)$ as the state, where each column and constraint is represented by node 
$v\in V$ and there is an edge $e \in E$ between columns and constraints only if the column contributes to the constraint. Detailed encoding for state is discussed in Section \ref{sec421}. As we are using a bipartite graph, the GNN we use should be able to achieve convolution on a bipartite graph with two types of nodes (variable nodes and constraint nodes), and we utilize the similar bipartite GNN as in the study conducted by (\cite{morabit}) with modification of the task it performs (from binary classification to Q value regression). Here we give a brief overview of how convolution or the features update is achieved in this bipartite GNN. 

The features update is done in two phase: phase 1 updates the constraint features and phase 2 updates the column features. In phase 1, the constraint features in next iteration is obtained by applying a non-linear transformation of previous constraint node features and aggregate information of its neighbouring nodes, which are column nodes that are connected to this constraint node. This can be treated as information passing from variable nodes to constraint nodes. Similarly,the second phase can be seen as message passing from constraint nodes to variable nodes. 

Once the column node features have been updated for several iterations for all the column nodes, these features are fed into a fully connected layer, which results in Q values for each column node.

\section{Training, Validation, Testing set}
\label{Append: Train,Val,Test set numbers}
Table \ref{tab:datasets CSP} below shows the information of the instances contained in the training set, validation set and testing set for CSP. Column lists total number of instances in each dataset, while other columns list the number of instances with specific roll length n in that dataset. The division of VRPTW dataset can be found in Section \ref{sec:VRPTW}.

{\renewcommand{\arraystretch}{1.1}
\begin{table}[htb]\small
\caption{Dataset division for \textbf{CSP}}
\label{tab:datasets CSP}
\centering
\begin{tabular}{|m{15mm}|
>{\centering}m{10mm}|
>{\centering}m{7mm}|
>{\centering}m{7mm}|
>{\centering}m{7mm}|
>{\centering\arraybackslash}m{7mm}|}\hline
Dataset & Total & $n=50$  &  $n=100 $  &  $n=200$   &  $n=750$  \\ \hline
Training & 400 & 160 & 160 & 80 & 0 \\
Validation & 30 & 10 & 10 & 10 & 0 \\
Testing & 156 & 49 & 0 & 86 & 21 \\\hline
\end{tabular}
\end{table}
}


\section{Curriculum Learning design}

In Table \ref{tab:learn}, we provide data characteristic we considered, for the sake of curriculum learning of the RL agent. In this paper the RL agent is trained on instances that are ordered according to their difficulty level. To accomplish this instances are divided into three categories of easy, normal and hard according to the stock length for CSP. 

Easy instances have stock length of 50, normal instances have stock length of 100 and hard instances have stock length of 200. There are 40 instances for each instance type. Details of training curriculum  of different instance types are shown in table \ref{tab:learn}. Figure \ref{fig:no_s_vs_s} displays the number of steps to convergence vs. instance number for training instances. It is clear that for each instance type the steps taken for convergence decreases as the model is trained on more of the same instance type. This shows that the RL agent successfully learns to select columns to enter basis. However, when instances are ordered randomly (Figure \ref{fig:b}) there are no specific trend on the steps taken to converge. This highlights the necessity of curriculum learning. Curriculum for VRPTW can be found in Section \ref{sec:VRPTW}.%

{\renewcommand{\arraystretch}{1.1}
\begin{table}[htb]\small
\caption{Curriculum Learning schedule for \textbf{CSP}}
\label{tab:learn}\centering
\begin{tabular}{|m{22mm}|
>{\centering}m{15mm}|
>{\centering}m{13mm}|
>{\centering\arraybackslash}m{13mm}|}\hline
\multicolumn{4}{|c|}{Training curriculum} 
\\\hline
Type of Instance & Number of Instances & Stock Length & 
Number of Orders \\\hline
Easy & 40 & 50 & 50\\
Easy & 40 & 50 & 75   \\
Easy & 40 & 50 & 100 \\
Easy & 40 & 50 & 120  \\
Normal & 40 & 100 & 75 \\
Normal & 40 & 100 & 100  \\
Normal & 40 & 100 & 120 \\
Normal & 40 & 100 & 150 \\
Hard & 40 & 200 & 125 \\
Hard &40 & 200 & 150 \\
\hline
\end{tabular}
\end{table}
}

\section{Node Features}
\label{Features}
Node features used for \textbf{CSP}:
\begin{enumerate}
    \item \textbf{Column node features}:\\
    Feature (a) and (c) relate to solving the RMP problem as they are all information about decision variables in RMP, and  each column node corresponds to one decision variable. Feature (b) and (d) are determined by the problem formulation of each cutting stock instance, while feature (e) - (i) corresponds to the dynamical information of each column entering and leaving the basis.
    \begin{enumerate}
    \item \textbf{Reduced cost}: Reduced cost is a quantity associated with each variable indicating how much the objective function coefficient on the corresponding variable must be improved before the solution value of the decision variable will be positive in the optimal solution (the cutting pattern will be used in optimal set of cutting patterns). The reduced cost value is only non-zero when the optimal value of a variable is zero.

    \item \textbf{Connectivity of column node}: Total number of constraint nodes each column node connects to. As each constraint is a particular demand, this node feature indicates the ability of a column node (a pattern) to satisfy demands. It also indicates the connectivity of each column node in the bipartite graph representing the state.
    
    \item \textbf{Solution value}: The solution value of each decision variable corresponding to each column node after solving the RMP in the current iteration.
    For each column node, this feature is continuous number greater than or equal to 0. The candidate column nodes have this feature set to be 0.
    \item \textbf{Waste}:  A feature recording the remaining length of a roll if we were to cut the current pattern from the roll. Again, each column node corresponds to one decision variable in RMP, which also represents one particular cut pattern.
    \item \textbf{Number of iterations that each column node stays in the basis}:  If the column node stays in the basis for a long time, it is most likely that the pattern corresponds to this column node is really good.
    \item \textbf{Number of iterations that each column node stays out of the basis}: if the column node stays out of the basis for a long time, it is most likely never enters the basis and being used in optimal solution in future iterations.  
    \item \textbf{If the column left basis on the last iteration or not}: This is a binary feature recording the dynamics of each column node.
    \item \textbf{If the column enter basis on the last iteration or not}: Similar binary feature as (f).
    \item \textbf{Action node or not}: A binary feature indicating whether a column node is a candidate (a newly added action) or not. If the column node is a candidate node (column) that is generated at the current iteration by SP, then this binary feature is 1 otherwise 0.
\end{enumerate}
    \item \textbf{Constraint node features}:\\
       Each constraint node corresponds to one constraint  in RMP, so the number of constraint nodes are fixed for each cutting stock problem instance.
\begin{enumerate}
\item \textbf{Dual value}: Dual value or shadow price is the coefficient of each dual variable in sub-problem objective function, and as each constraint node corresponds to one dual variable, we record dual value as one feature for constraint node.

\item \textbf{Connectivity of constraint node}: Total number of column nodes each constraint node connects to, which also indicates the connectivity of each constraint node in the bipartite graph representing the state.
\end{enumerate}
\end{enumerate}

\textbf{VRPTW} node features used are quite similar to CSP:\\ 

\begin{enumerate}
\item \textbf{Column node features}: Reduced cost, connectivity of column node, solution value, route cost, Number of iterations that each column node stays in the basis, Number of iterations that each column node stays in the basis, If the column left basis on the last iteration or no, If the column enter basis on the last iteration or no.

\item\textbf{Constraint node features}: Dual value, connectivity of constraint node
\end{enumerate}

\section{Detailed statistics of testing results}
\label{app:testing_results_more}
In Table \ref{tab:statistics for CSP}, we provide statistics obtained from our experimental results for CSP. We report both the average and standard deviation of number of iterations, solution time measured in seconds, and the objective function values. Note that the resulting objective function values between column selection policies might differ due to the early-stopping criteria we adopted. Note that for the sake table spacing, we use Objval to denote objective function value, and $\mu, \sigma$ for mean and standard deviation, respectively. We observe the clear dominance of~\method{} over the potential of~\method{} in solving challenging CG problems in practice.

{\tabcolsep=0.28\tabcolsep
\renewcommand{\arraystretch}{1.3}
\begin{table*}[htb]\footnotesize
\centering
\begin{tabular}{|c|cccccc|cccccc|cccccc|cccccc|}
\hline
& \multicolumn{6}{c|}{$n=50$}                                                           
& \multicolumn{6}{c|}{$n=200$}                                                                       &  \multicolumn{6}{c|}{$n=750$}                                                                                           \\ \hline
                & \multicolumn{2}{c|}{\textbf{Iteration}} & \multicolumn{2}{c|}{\textbf{Time(s)}} & \multicolumn{2}{c|}{\textbf{Objval}} & \multicolumn{2}{c|}{\textbf{Iteration}} & \multicolumn{2}{c|}{\textbf{Time(s)}} & \multicolumn{2}{c|}{\textbf{Objval}}  & \multicolumn{2}{c|}{\textbf{Iteration}} & \multicolumn{2}{c|}{\textbf{Time(s)}} & \multicolumn{2}{c|}{\textbf{Objval}} \\ \hline
                & $\mu$  & \multicolumn{1}{c|}{$\sigma$}  & $\mu$ & \multicolumn{1}{c|}{$\sigma$} & $\mu$           & $\sigma$           & $\mu$  & \multicolumn{1}{c|}{$\sigma$}  & $\mu$ & \multicolumn{1}{c|}{$\sigma$} & $\mu$           & $\sigma$           & $\mu$  & \multicolumn{1}{c|}{$\sigma$}  & $\mu$ & \multicolumn{1}{c|}{$\sigma$} & $\mu$           & $\sigma$         
                \\ \hline
\textbf{Greedy} & 51.2   & \multicolumn{1}{c|}{11.4}      & 6.0   & \multicolumn{1}{c|}{2.5}      & 23.4            & 2.8                & 78.0   & \multicolumn{1}{c|}{23.0}      & 22.3  & \multicolumn{1}{c|}{14.4}     & 91.1            & 10.1                            & 292.8  & \multicolumn{1}{c|}{47.7}      & 640.8 & \multicolumn{1}{c|}{240.2}     & 327.6           & 18.2              \\
\textbf{RL}     & 41.6   & \multicolumn{1}{c|}{10.0}       & 4.9   & \multicolumn{1}{c|}{2.2}      & 23.5            & 2.8                & 64.7   & \multicolumn{1}{c|}{20.5}      & 18.1  & \multicolumn{1}{c|}{12.3}      & 91.5            & 9.8                          & 227.2  & \multicolumn{1}{c|}{42.7}      & 460.3 & \multicolumn{1}{c|}{219.9}     & 328.6           & 18.2   
\\
\textbf{Expert} & 42.5   & \multicolumn{1}{c|}{9.5}      & 7.1   & \multicolumn{1}{c|}{3.1}      & 23.4            & 2.8                & 69.1   & \multicolumn{1}{c|}{18.8}      & 29.3  & \multicolumn{1}{c|}{18.9}     & 91.2            &  10.1          & 251.3  & \multicolumn{1}{c|}{35.2}      & 827.9 & \multicolumn{1}{c|}{302.0}     & 327.2           & 18.7  \\\hline

\end{tabular}
\caption{Solution time, iteration and objective function value reports with $\mu$ mean, and $\sigma$ standard deviation for \text{CSP}}
\label{tab:statistics for CSP}
\end{table*}}

In Table \ref{tab:statistics for VRP}, we provide the same statistics obtained from our experimental results for VRPTW.
{\tabcolsep=0.28\tabcolsep
\renewcommand{\arraystretch}{1.3}
\begin{table*}[htb]\footnotesize
\centering

\begin{tabular}{|c|cccccc|cccccc|cccccc|}
\hline
& \multicolumn{6}{c|}{small test instances}                                                           
& \multicolumn{6}{c|}{large test instances}                                                                              \\ \hline
                & \multicolumn{2}{c|}{\textbf{Iteration}} & \multicolumn{2}{c|}{\textbf{Time(s)}} & \multicolumn{2}{c|}{\textbf{Objval}} & \multicolumn{2}{c|}{\textbf{Iteration}} & \multicolumn{2}{c|}{\textbf{Time(s)}} & \multicolumn{2}{c|}{\textbf{Objval}} \\ \hline
                & $\mu$  & \multicolumn{1}{c|}{$\sigma$}  & $\mu$ & \multicolumn{1}{c|}{$\sigma$} & $\mu$           & $\sigma$           & $\mu$  & \multicolumn{1}{c|}{$\sigma$}  & $\mu$ & \multicolumn{1}{c|}{$\sigma$} & $\mu$           & $\sigma$       
                \\ \hline
\textbf{Greedy} & 18.8   & \multicolumn{1}{c|}{ 11.7}      &  291.1   & \multicolumn{1}{c|}{ 285.0}      &  102.6  &  30 &  128.1   & \multicolumn{1}{c|}{ 100.3}      &  3268.2  & \multicolumn{1}{c|}{ 4627.5}     &  458.6            &  144.6                             \\
\textbf{RL}     &  9.5   & \multicolumn{1}{c|}{6.3}       &  179.9   & \multicolumn{1}{c|}{ 247.1}      &  102.6            &  30.0                &  75.6   & \multicolumn{1}{c|}{ 39.3}      &  1832.8  & \multicolumn{1}{c|}{ 2015.3}      &  458.6           &  144.6                             

\\\hline

\end{tabular}
\caption{Solution time, iteration and objective function value reports with $\mu$ mean, and $\sigma$ standard deviation for \textbf{VRPTW}}
\label{tab:statistics for VRP}
\end{table*}}

\section{Per instance results CSP}
\label{app:per_ins_result}

\begin{tabular}{lrrrrrr}
\toprule
             Names & \multicolumn{2}{l}{Iterations} & \multicolumn{2}{l}{Time} & \multicolumn{2}{l}{Objective value} \\
     Instance name & Greedy\_iter & RL\_iter & Greedy\_time & RL\_time &      Greedy\_obj & RL\_obj \\
\midrule
 r107.txt,  n = 20 &         328 &      84 &    15164.16 & 3843.84 &          517.00 & 517.00 \\
 r110.txt,  n = 19 &         189 &     101 &     4093.69 & 2170.76 &          522.00 & 522.00 \\
 c101.txt,  n = 25 &         140 &      90 &     4123.05 & 2632.07 &          637.00 & 637.00 \\
 r106.txt,  n = 20 &          56 &      30 &      807.00 &  429.26 &          637.00 & 637.00 \\
 r111.txt,  n = 20 &         203 &     126 &     6965.91 & 4350.76 &          534.00 & 534.00 \\
 r104.txt,  n = 20 &         261 &     119 &    18788.82 & 8359.52 &          504.00 & 504.00 \\
 r108.txt,  n = 18 &          49 &      39 &     6062.88 & 4529.36 &          156.00 & 156.00 \\
 r103.txt,  n = 19 &         440 &     142 &    15789.96 & 5132.63 &          437.00 & 437.00 \\
rc104.txt,  n = 15 &         196 &      79 &     2992.44 & 1211.46 &          491.00 & 491.00 \\
rc103.txt,  n = 19 &         163 &      84 &     2135.99 & 1106.87 &          667.00 & 667.00 \\
rc107.txt,  n = 15 &          77 &      77 &      766.64 &  783.29 &          341.00 & 341.00 \\
rc105.txt,  n = 15 &         107 &      88 &      578.83 &  497.95 &          415.00 & 415.00 \\
rc101.txt,  n = 16 &          40 &      27 &      106.21 &   69.52 &          571.00 & 571.00 \\
rc102.txt,  n = 14 &         115 &      73 &      999.08 &  645.00 &          391.00 & 391.00 \\
rc108.txt,  n = 15 &         345 &     173 &     5449.42 & 2649.93 &          383.00 & 383.00 \\
 r102.txt,  n = 14 &         110 &      59 &      960.38 &  513.88 &          382.00 & 382.00 \\
 c101.txt,  n = 15 &          36 &      32 &       88.32 &   78.41 &          491.00 & 491.00 \\
 r105.txt,  n = 15 &         173 &      93 &     3877.94 & 2100.79 &          583.00 & 583.00 \\
 c107.txt,  n = 15 &          93 &      83 &     5609.29 & 4938.09 &          249.00 & 249.00 \\
 c105.txt,  n = 20 &          87 &      81 &     3874.88 & 3609.85 &          485.00 & 485.00 \\
rc201.txt,  n = 16 &          57 &      51 &     1752.45 & 1567.67 &          641.00 & 641.00 \\
rc106.txt,  n = 15 &         128 &      98 &      991.33 &  745.47 &          276.00 & 276.00 \\
 r101.txt,  n = 15 &          63 &      55 &      100.00 &   87.98 &          409.00 & 409.00 \\
 r205.txt,  n = 19 &          11 &      16 &     1498.44 & 2213.35 &          749.00 & 749.00 \\
 r108.txt,  n = 15 &          37 &      30 &     2992.61 & 2321.59 &          148.75 & 148.75 \\
 r109.txt,  n = 19 &         153 &      71 &     1294.90 &  603.83 &          578.00 & 578.00 \\

\bottomrule
\end{tabular}

\section{VRPTW box plots}
\label{sec:VRPTW other plots}
Similar to CSP, we see the upward generalization of~\method{} in VRPTW: we only train our model with small sized instances, and tested on large instances with more customers, and we observe that~\method{} was able to perform well using VRPTW testing instances with large size. These results, again indicate that our proposed~\method{} can be preferable in solving challenging column generation problems because of its ability to generalize well. Besides, compare large testing results with small testing results for VRPTW (also compare n=750 results with n=50,200 for CSP), we observe that the more challenging the problem,  the larger gap exist between~\method{} and benchmark method. This indicates that for CG problem we considered, the harder the problem, the more important 
considering the future effect of CG becomes, thus taking future into consideration would drastically accelerate the solving process of CG.
\begin{figure}[!htp]
\centering
\includegraphics[width=.30\textwidth]{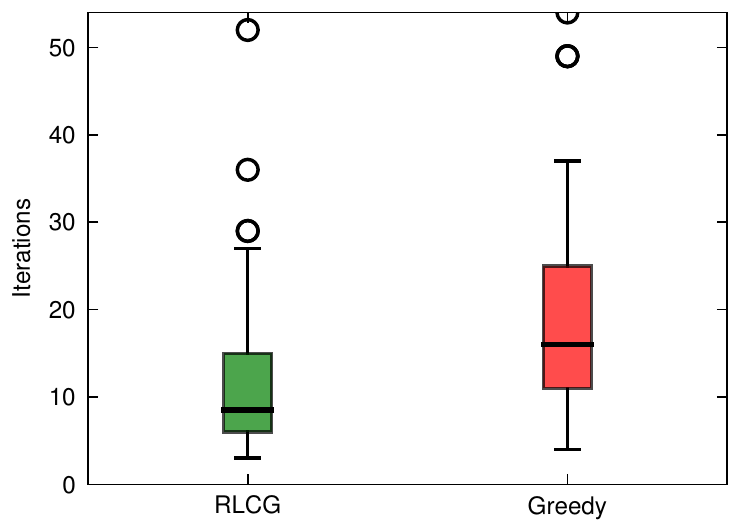}\quad
\includegraphics[width=.30\textwidth]{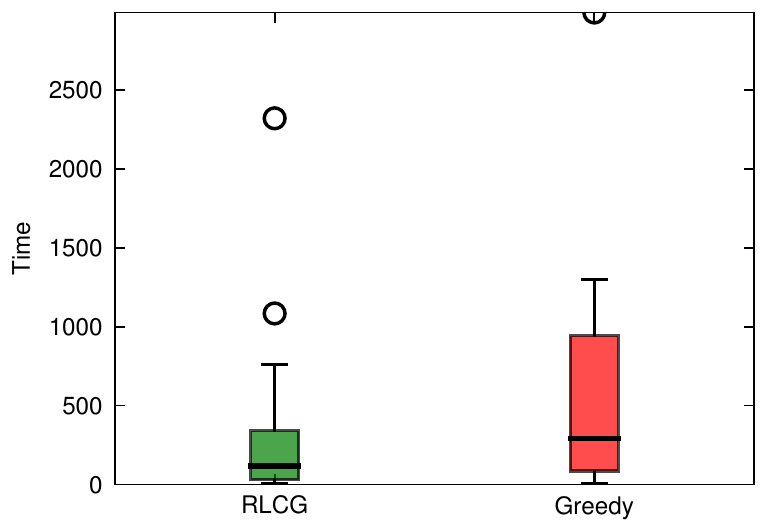} 
\caption{Pair-wise comparisons between~\method{} and greedy benchmark in terms of the number of CG iterations and the solving time over small and large testing sets for \textbf{VRPTW}}
\label{fig:box-VRP-small}
\end{figure}

\section{per instance results CSP}
\label{app:per_ins_result_CSP}
\begin{table*}[t]
\centering
\resizebox{\textwidth}{!}{
\begin{tabular}{l|ccc|rrr|ccc}
\toprule
Name & \multicolumn{3}{l}{Iterations} & \multicolumn{3}{l}{Time} & \multicolumn{3}{l}{Objective value} \\
Instance name & Greedy\_iter & RL\_iter & Expert\_iter & Greedy\_time & RL\_time & Expert\_time &      Greedy\_obj & RL\_obj & Expert\_obj \\
\midrule
  BPP\_50\_125\_0.1\_0.7\_2 &          45 &      30 &          28 &        4.45 &    2.80 &        3.36 &           18.18 &  18.32 &      18.31 \\
  BPP\_50\_200\_0.2\_0.8\_5 &          46 &      31 &          31 &        5.15 &    3.46 &        5.12 &           29.00 &  29.00 &      29.00 \\
  BPP\_50\_50\_0.1\_0.8\_0 &          40 &      31 &          31 &        2.64 &    2.22 &        2.81 &           23.50 &  23.50 &      23.50 \\
  BPP\_50\_125\_0.1\_0.8\_8 &          53 &      38 &          48 &        6.04 &    4.31 &        7.84 &           20.74 &  21.31 &      20.64 \\
  BPP\_50\_200\_0.2\_0.7\_2 &          57 &      42 &          38 &        6.70 &    5.10 &        6.47 &           23.06 &  23.00 &      23.00 \\
  BPP\_50\_125\_0.1\_0.7\_8 &          57 &      52 &          47 &        6.13 &    5.65 &        6.82 &           20.04 &  20.31 &      20.18 \\
  BPP\_50\_200\_0.2\_0.7\_0 &          43 &      42 &          35 &        4.48 &    4.85 &        5.50 &           23.46 &  23.58 &      23.45 \\
  BPP\_50\_125\_0.1\_0.8\_6 &          61 &      43 &          57 &        8.35 &    5.71 &       11.61 &           25.25 &  25.50 &      25.25 \\
  BPP\_50\_200\_0.2\_0.8\_7 &          57 &      44 &          50 &        6.08 &    4.83 &        7.93 &           24.00 &  23.96 &      23.94 \\
  BPP\_50\_125\_0.1\_0.8\_1 &          67 &      54 &          49 &        7.64 &    6.36 &        7.83 &           22.67 &  22.79 &      22.68 \\
  BPP\_50\_50\_0.1\_0.8\_1 &          25 &      23 &          23 &        1.50 &    1.56 &        1.92 &           28.00 &  28.00 &      28.00 \\
  BPP\_50\_200\_0.2\_0.7\_3 &          63 &      45 &          54 &        7.91 &    5.44 &        9.92 &           23.50 &  23.50 &      23.50 \\
  BPP\_50\_125\_0.2\_0.7\_4 &          44 &      37 &          31 &        3.08 &    2.80 &        3.02 &           25.00 &  25.00 &      25.00 \\
  BPP\_50\_200\_0.2\_0.7\_5 &          68 &      52 &          46 &        9.80 &    7.33 &        9.08 &           24.60 &  24.70 &      24.70 \\
  BPP\_50\_125\_0.1\_0.8\_3 &          35 &      29 &          34 &        3.42 &    2.90 &        4.97 &           28.50 &  28.50 &      28.50 \\
  BPP\_50\_200\_0.1\_0.8\_2 &          59 &      55 &          52 &        9.66 &    9.25 &       12.15 &           20.29 &  20.51 &      20.33 \\
  BPP\_50\_125\_0.1\_0.7\_1 &          46 &      48 &          37 &        5.11 &    5.79 &        5.79 &           18.92 &  18.95 &      18.98 \\
  BPP\_50\_200\_0.2\_0.7\_8 &          56 &      31 &          34 &        7.96 &    3.94 &        6.69 &           21.65 &  21.83 &      21.62 \\
  BPP\_50\_200\_0.2\_0.7\_9 &          50 &      38 &          42 &        5.37 &    4.00 &        6.39 &           23.21 &  23.56 &      23.15 \\
  BPP\_50\_200\_0.1\_0.8\_1 &          62 &      56 &          59 &        8.62 &    8.05 &       12.62 &           23.21 &  23.25 &      23.07 \\
  BPP\_50\_200\_0.1\_0.8\_4 &          58 &      52 &          46 &        8.05 &    7.25 &        9.00 &           25.00 &  25.00 &      25.00 \\
  BPP\_50\_125\_0.1\_0.7\_4 &          51 &      42 &          41 &        6.13 &    5.01 &        6.78 &           21.30 &  21.33 &      21.42 \\
  BPP\_50\_125\_0.2\_0.7\_0 &          44 &      41 &          46 &        3.90 &    3.91 &        6.26 &           21.38 &  21.62 &      21.38 \\
  BPP\_50\_200\_0.1\_0.8\_6 &          62 &      50 &          53 &        8.99 &    6.96 &       10.77 &           22.60 &  22.86 &      22.58 \\
  BPP\_50\_200\_0.1\_0.8\_0 &          84 &      69 &          65 &       14.12 &   11.63 &       14.99 &           23.17 &  23.57 &      23.15 \\
  BPP\_50\_125\_0.1\_0.7\_5 &          45 &      28 &          36 &        4.90 &    2.87 &        5.43 &           19.71 &  20.04 &      19.88 \\
  BPP\_50\_200\_0.2\_0.8\_2 &          53 &      48 &          47 &        6.58 &    6.21 &        8.93 &           26.00 &  26.00 &      26.00 \\
  BPP\_50\_50\_0.1\_0.8\_2 &          31 &      29 &          30 &        1.62 &    1.66 &        2.06 &           25.00 &  25.00 &      25.00 \\
  BPP\_50\_125\_0.2\_0.7\_2 &          46 &      38 &          39 &        3.98 &    3.29 &        4.66 &           21.74 &  22.07 &      21.89 \\
  BPP\_50\_50\_0.1\_0.8\_4 &          39 &      35 &          34 &        2.56 &    2.53 &        3.21 &           23.25 &  23.29 &      23.25 \\
  BPP\_50\_200\_0.2\_0.7\_4 &          55 &      33 &          39 &        5.51 &    3.17 &        5.24 &           21.83 &  21.92 &      21.83 \\
  BPP\_50\_125\_0.1\_0.8\_7 &          50 &      43 &          46 &        5.70 &    5.20 &        8.07 &           26.00 &  26.00 &      26.00 \\
  BPP\_50\_200\_0.1\_0.8\_7 &          59 &      75 &          66 &        8.85 &   12.93 &       15.27 &           21.36 &  21.21 &      21.19 \\
  BPP\_50\_125\_0.2\_0.7\_3 &          46 &      39 &          38 &        3.52 &    3.19 &        3.95 &           23.00 &  23.00 &      23.00 \\
  BPP\_50\_125\_0.1\_0.8\_2 &          39 &      40 &          42 &        4.14 &    4.73 &        7.02 &           20.80 &  20.82 &      20.62 \\
  BPP\_50\_125\_0.1\_0.7\_7 &          47 &      42 &          41 &        5.65 &    5.04 &        6.99 &           19.63 &  19.81 &      19.70 \\
  BPP\_50\_125\_0.1\_0.8\_5 &          40 &      36 &          35 &        3.67 &    3.46 &        4.66 &           29.50 &  29.50 &      29.50 \\
  BPP\_50\_125\_0.2\_0.7\_1 &          61 &      44 &          48 &        7.03 &    4.97 &        7.88 &           23.15 &  23.15 &      23.15 \\
  BPP\_50\_200\_0.2\_0.8\_0 &          70 &      42 &          51 &       10.33 &    5.75 &       10.72 &           27.50 &  27.50 &      27.50 \\
  BPP\_50\_125\_0.1\_0.8\_9 &          50 &      43 &          48 &        5.94 &    5.35 &        8.84 &           23.50 &  23.50 &      23.50 \\
  BPP\_50\_125\_0.1\_0.7\_6 &          52 &      40 &          40 &        6.01 &    4.64 &        6.51 &           20.77 &  20.98 &      20.90 \\
  BPP\_50\_200\_0.1\_0.8\_9 &          64 &      42 &          52 &        8.31 &    5.28 &        9.66 &           21.88 &  21.99 &      21.83 \\
  BPP\_50\_50\_0.1\_0.8\_3 &          30 &      30 &          26 &        1.65 &    1.76 &        1.75 &           20.06 &  20.12 &      20.07 \\
  BPP\_50\_125\_0.1\_0.7\_9 &          51 &      42 &          43 &        6.26 &    5.10 &        7.23 &           20.46 &  20.66 &      20.59 \\
  BPP\_50\_125\_0.1\_0.8\_0 &          49 &      41 &          46 &        5.32 &    4.70 &        7.57 &           26.00 &  26.00 &      26.00 \\
  BPP\_50\_200\_0.1\_0.8\_3 &          66 &      45 &          42 &        8.79 &    5.84 &        7.68 &           21.98 &  22.14 &      22.21 \\
  BPP\_50\_200\_0.1\_0.8\_5 &          34 &      33 &          35 &        3.48 &    3.57 &        5.60 &           30.00 &  30.00 &      30.00 \\
  BPP\_50\_200\_0.2\_0.7\_7 &          57 &      48 &          50 &        6.51 &    5.61 &        8.30 &           23.50 &  23.50 &      23.50 \\
  BPP\_50\_200\_0.2\_0.8\_1 &          41 &      28 &          33 &        4.82 &    3.32 &        6.06 &           27.50 &  27.50 &      27.50 \\

\bottomrule
\end{tabular}}
\end{table*}

\clearpage
\begin{table*}[t]
\centering
\resizebox{\textwidth}{!}{
\begin{tabular}{l|rrr|rrr|rrr}
\toprule
Name & \multicolumn{3}{l}{Iterations} & \multicolumn{3}{l}{Time} & \multicolumn{3}{l}{Objective value} \\
Instance name & Greedy\_iter & RL\_iter & Expert\_iter & Greedy\_time & RL\_time & Expert\_time &      Greedy\_obj & RL\_obj & Expert\_obj \\
\midrule
 BPP\_200\_100\_0.2\_0.7\_1 &          57 &      45 &          56 &       11.86 &    9.18 &       18.19 &           90.50 &  90.81 &      90.50 \\
 BPP\_200\_100\_0.1\_0.7\_4 &          80 &      61 &          68 &       24.36 &   18.23 &       31.03 &           79.53 &  79.93 &      79.70 \\
  BPP\_200\_75\_0.2\_0.8\_5 &          63 &      52 &          52 &       10.91 &    9.10 &       13.37 &          103.00 & 103.00 &     103.00 \\
  BPP\_200\_75\_0.2\_0.7\_9 &          55 &      47 &          52 &        7.55 &    6.43 &       10.35 &           89.24 &  89.67 &      89.25 \\
  BPP\_200\_75\_0.1\_0.7\_4 &          66 &      48 &          67 &       12.67 &    8.93 &       19.91 &           80.44 &  81.72 &      80.37 \\
 BPP\_200\_100\_0.2\_0.7\_0 &          65 &      53 &          59 &       13.45 &   11.01 &       19.05 &           93.45 &  93.40 &      93.40 \\
 BPP\_200\_100\_0.1\_0.7\_5 &          69 &      49 &          71 &       19.29 &   13.06 &       31.68 &           78.27 &  79.52 &      78.02 \\
  BPP\_200\_50\_0.2\_0.8\_1 &          30 &      30 &          29 &        2.42 &    2.70 &        3.44 &          115.50 & 115.50 &     115.50 \\
 BPP\_200\_120\_0.1\_0.7\_7 &          87 &      72 &          87 &       33.18 &   26.98 &       52.32 &           78.47 &  79.01 &      78.47 \\
  BPP\_200\_50\_0.2\_0.8\_5 &          44 &      32 &          42 &        3.92 &    2.93 &        5.46 &          103.00 & 103.00 &     103.00 \\
 BPP\_200\_100\_0.1\_0.7\_7 &          71 &      56 &          65 &       20.23 &   15.24 &       28.03 &           77.10 &  77.60 &      77.03 \\
  BPP\_200\_75\_0.1\_0.7\_7 &          78 &      57 &          64 &       16.39 &   11.12 &       18.52 &           82.80 &  83.70 &      82.92 \\
  BPP\_200\_75\_0.1\_0.7\_1 &          75 &      63 &          69 &       15.37 &   12.84 &       21.38 &           82.69 &  83.60 &      82.64 \\
  BPP\_200\_75\_0.1\_0.8\_7 &          97 &      82 &          82 &       26.72 &   21.65 &       31.89 &           91.71 &  92.08 &      91.72 \\
  BPP\_200\_75\_0.2\_0.8\_0 &          59 &      48 &          59 &        9.94 &    8.15 &       15.58 &          102.33 & 102.33 &     102.33 \\
 BPP\_200\_120\_0.1\_0.7\_3 &         111 &      89 &          92 &       51.56 &   39.32 &       63.39 &           82.81 &  83.30 &      83.24 \\
 BPP\_200\_100\_0.2\_0.8\_8 &          85 &      67 &          73 &       24.35 &   18.42 &       31.85 &          111.00 & 111.00 &     111.00 \\
 BPP\_200\_100\_0.2\_0.7\_8 &          88 &      51 &          75 &       20.32 &   10.12 &       25.20 &           91.21 &  92.33 &      91.19 \\
 BPP\_200\_100\_0.2\_0.8\_2 &         107 &      89 &          88 &       33.04 &   26.74 &       39.09 &           99.56 &  99.59 &      99.56 \\
  BPP\_200\_75\_0.2\_0.7\_2 &          56 &      46 &          49 &        7.75 &    6.12 &        9.31 &           92.30 &  92.96 &      92.44 \\
  BPP\_200\_75\_0.1\_0.8\_9 &          84 &      87 &          87 &       20.47 &   23.16 &       33.23 &           92.00 &  92.28 &      92.00 \\
 BPP\_200\_100\_0.1\_0.8\_4 &         106 &      87 &          85 &       43.36 &   34.82 &       51.75 &           94.50 &  94.50 &      94.50 \\
 BPP\_200\_100\_0.1\_0.8\_6 &         105 &      84 &          87 &       39.42 &   30.47 &       47.84 &           88.96 &  89.42 &      89.18 \\
  BPP\_200\_75\_0.1\_0.8\_4 &          84 &      82 &          75 &       21.36 &   21.97 &       28.12 &           92.03 &  92.06 &      92.01 \\
 BPP\_200\_100\_0.2\_0.8\_9 &         121 &      92 &          89 &       40.25 &   28.35 &       40.78 &          101.94 & 102.00 &     101.87 \\
  BPP\_200\_75\_0.2\_0.8\_8 &          75 &      55 &          59 &       13.82 &    9.63 &       15.57 &          101.50 & 101.50 &     101.50 \\
 BPP\_200\_100\_0.1\_0.8\_8 &         124 &     112 &         106 &       56.95 &   52.58 &       73.34 &           93.24 &  93.29 &      93.58 \\
  BPP\_200\_75\_0.1\_0.7\_2 &          60 &      48 &          54 &       11.35 &    8.97 &       15.03 &           78.03 &  78.43 &      78.03 \\
  BPP\_200\_75\_0.2\_0.8\_2 &          71 &      55 &          61 &       12.99 &    9.86 &       16.48 &           96.65 &  96.58 &      96.58 \\
  BPP\_200\_75\_0.2\_0.7\_7 &          54 &      50 &          50 &        7.37 &    6.85 &        9.84 &           91.00 &  91.00 &      91.00 \\
 BPP\_200\_100\_0.1\_0.7\_9 &          77 &      70 &          77 &       24.27 &   22.05 &       37.97 &           80.15 &  80.62 &      80.34 \\
  BPP\_200\_75\_0.1\_0.7\_9 &          58 &      54 &          65 &       10.46 &    9.92 &       18.38 &           82.59 &  83.25 &      82.41 \\
 BPP\_200\_100\_0.2\_0.8\_6 &         110 &      93 &          99 &       36.26 &   29.91 &       48.23 &           98.00 &  98.00 &      98.00 \\
  BPP\_200\_75\_0.2\_0.8\_3 &          61 &      50 &          58 &       10.39 &    8.45 &       15.06 &          102.50 & 102.50 &     102.50 \\
 BPP\_200\_100\_0.1\_0.7\_2 &          95 &      69 &          72 &       31.75 &   20.94 &       33.39 &           80.63 &  81.30 &      81.02 \\
  BPP\_200\_75\_0.2\_0.7\_8 &          56 &      39 &          50 &        7.73 &    4.96 &        9.78 &           89.12 &  89.70 &      89.12 \\
 BPP\_200\_100\_0.1\_0.7\_8 &          85 &      74 &          85 &       27.10 &   23.43 &       41.74 &           80.91 &  81.01 &      80.70 \\
 BPP\_200\_100\_0.1\_0.8\_5 &         102 &      86 &          95 &       42.48 &   34.58 &       61.02 &           90.08 &  91.64 &      89.99 \\
 BPP\_200\_100\_0.1\_0.7\_3 &          80 &      67 &          65 &       22.67 &   18.24 &       26.34 &           79.88 &  80.18 &      79.66 \\
  BPP\_200\_75\_0.2\_0.8\_6 &          54 &      41 &          43 &        8.81 &    6.63 &       10.44 &          115.00 & 115.00 &     115.00 \\
  BPP\_200\_75\_0.2\_0.8\_7 &          57 &      52 &          54 &        9.08 &    8.60 &       13.34 &          108.00 & 108.00 &     108.00 \\
 BPP\_200\_120\_0.1\_0.7\_0 &         104 &      84 &          81 &       44.38 &   34.28 &       50.55 &           80.17 &  80.86 &      80.89 \\
  BPP\_200\_75\_0.1\_0.7\_4 &          66 &      48 &          67 &       12.58 &    8.78 &       19.99 &           80.44 &  81.72 &      80.37 \\
 BPP\_200\_100\_0.2\_0.7\_7 &          87 &      72 &          78 &       21.06 &   16.82 &       27.69 &           93.57 &  93.69 &      93.50 \\
 BPP\_200\_120\_0.1\_0.7\_8 &          69 &      53 &          73 &       22.31 &   16.40 &       38.58 &           79.73 &  80.62 &      79.57 \\
 BPP\_200\_100\_0.1\_0.8\_3 &         106 &      93 &          88 &       40.46 &   35.01 &       49.13 &           89.03 &  89.98 &      89.15 \\
  BPP\_200\_75\_0.1\_0.8\_8 &          87 &      81 &          81 &       23.46 &   22.07 &       32.42 &           87.02 &  87.65 &      87.03 \\
 BPP\_200\_100\_0.1\_0.7\_0 &          80 &      72 &          74 &       23.89 &   21.35 &       33.65 &           81.23 &  81.34 &      80.98 \\
  BPP\_200\_50\_0.2\_0.8\_2 &          35 &      29 &          36 &        3.00 &    3.35 &        4.57 &          103.50 & 103.50 &     103.50 \\
 BPP\_200\_120\_0.1\_0.8\_1 &          99 &     102 &         102 &       48.69 &   52.01 &       81.43 &           84.03 &  84.11 &      83.68 \\
 BPP\_200\_120\_0.1\_0.7\_6 &          88 &      75 &          80 &       31.71 &   26.87 &       44.49 &           82.98 &  83.25 &      83.06 \\
 BPP\_200\_120\_0.1\_0.7\_5 &          95 &      78 &          78 &       37.97 &   30.99 &       46.90 &           80.48 &  80.68 &      80.82 \\
 BPP\_200\_100\_0.2\_0.8\_3 &          97 &      77 &          84 &       29.38 &   22.22 &       38.00 &          100.33 & 100.33 &     100.33 \\
  BPP\_200\_50\_0.2\_0.8\_3 &          40 &      33 &          34 &        3.54 &    3.04 &        4.16 &          101.83 & 101.83 &     101.83 \\
 BPP\_200\_100\_0.2\_0.8\_1 &          87 &      59 &          73 &       25.16 &   15.68 &       31.66 &          106.00 & 106.00 &     106.00 \\
 BPP\_200\_120\_0.1\_0.7\_9 &          94 &      80 &          76 &       36.10 &   30.36 &       51.57 &           79.99 &  80.27 &      79.99 \\
\bottomrule
\end{tabular}}
\end{table*}

\clearpage
\begin{table*}[t]
\centering
\resizebox{\textwidth}{!}{
\begin{tabular}{l|rrr|rrr|rrr}
\toprule
Name & \multicolumn{3}{l}{Iterations} & \multicolumn{3}{l}{Time} & \multicolumn{3}{l}{Objective value} \\
Instance name & Greedy\_iter & RL\_iter & Expert\_iter & Greedy\_time & RL\_time & Expert\_time &      Greedy\_obj & RL\_obj & Expert\_obj \\
\midrule
  BPP\_200\_75\_0.2\_0.7\_6 &          52 &      45 &          44 &        6.17 &    5.82 &        7.85 &           94.45 &  94.45 &      94.45 \\
  BPP\_200\_75\_0.1\_0.8\_5 &          99 &      84 &          74 &       27.22 &   22.59 &       28.18 &           89.04 &  89.48 &      89.34 \\
  BPP\_200\_50\_0.2\_0.8\_0 &          41 &      32 &          37 &        3.57 &    2.90 &        4.70 &           98.33 &  98.33 &      98.33 \\
  BPP\_200\_50\_0.2\_0.8\_9 &          36 &      28 &          33 &        3.06 &    2.44 &        4.16 &          108.50 & 108.50 &     108.50 \\
 BPP\_200\_120\_0.1\_0.7\_2 &         107 &      82 &          91 &       46.05 &   32.88 &       58.13 &           80.49 &  81.01 &      80.66 \\
  BPP\_200\_50\_0.2\_0.8\_6 &          37 &      30 &          36 &        3.23 &    2.64 &        4.54 &          101.50 & 101.50 &     101.50 \\
  BPP\_200\_50\_0.2\_0.8\_4 &          38 &      31 &          34 &        3.30 &    2.88 &        4.24 &          105.50 & 105.50 &     105.50 \\
  BPP\_200\_50\_0.2\_0.8\_7 &          36 &      35 &          32 &        3.00 &    3.23 &        3.90 &          102.00 & 102.00 &     102.00 \\
  BPP\_200\_75\_0.1\_0.8\_2 &          96 &      80 &          82 &       25.07 &   20.13 &       30.15 &           87.01 &  87.31 &      87.02 \\
 BPP\_200\_100\_0.1\_0.7\_6 &          97 &      72 &          81 &       31.30 &   21.39 &       37.84 &           81.77 &  82.77 &      82.23 \\
 BPP\_200\_100\_0.1\_0.8\_0 &         126 &     111 &         109 &       56.03 &   48.24 &       70.80 &           92.00 &  92.67 &      92.17 \\
  BPP\_200\_75\_0.1\_0.7\_0 &          73 &      59 &          68 &       14.63 &   11.65 &       20.31 &           83.04 &  83.35 &      82.97 \\
 BPP\_200\_100\_0.2\_0.8\_4 &          78 &      66 &          70 &       21.59 &   18.17 &       30.14 &          112.00 & 112.00 &     112.00 \\
  BPP\_200\_75\_0.2\_0.8\_4 &          87 &      65 &          75 &       17.02 &   12.03 &       20.98 &          102.73 & 102.73 &     102.73 \\
  BPP\_200\_75\_0.2\_0.7\_0 &          70 &      59 &          61 &       10.48 &    8.86 &       12.56 &           89.60 &  90.05 &      89.65 \\
  BPP\_200\_75\_0.1\_0.8\_6 &          81 &      72 &          74 &       20.60 &   17.78 &       27.62 &           84.92 &  85.28 &      84.83 \\
  BPP\_200\_75\_0.1\_0.8\_0 &          89 &      80 &          83 &       24.40 &   21.46 &       32.58 &           88.28 &  88.92 &      88.43 \\
  BPP\_200\_75\_0.1\_0.7\_8 &          56 &      52 &          64 &       10.33 &    9.84 &       19.50 &           77.35 &  77.48 &      77.21 \\
  BPP\_200\_75\_0.2\_0.8\_1 &          87 &      74 &          72 &       17.77 &   14.93 &       20.73 &          101.83 & 101.83 &     101.92 \\
  BPP\_200\_75\_0.1\_0.7\_6 &          74 &      61 &          57 &       14.43 &   12.27 &       16.06 &           82.75 &  82.68 &      82.81 \\
  BPP\_200\_75\_0.1\_0.8\_3 &          59 &      59 &          57 &       13.11 &   13.84 &       19.95 &          105.50 & 105.50 &     105.50 \\
 BPP\_200\_120\_0.1\_0.7\_1 &         104 &      86 &          84 &       43.51 &   34.34 &       51.44 &           82.28 &  82.78 &      82.28 \\
  BPP\_200\_75\_0.2\_0.8\_9 &          67 &      64 &          63 &       21.74 &   12.43 &       17.49 &           94.89 &  94.92 &      94.89 \\
 BPP\_200\_120\_0.1\_0.8\_3 &         118 &     100 &         104 &       61.68 &   50.27 &       81.97 &           89.24 &  89.80 &      89.14 \\
 BPP\_200\_120\_0.1\_0.7\_4 &          96 &      60 &          73 &       41.30 &   22.68 &       45.75 &           78.44 &  79.64 &      78.58 \\
 BPP\_200\_100\_0.2\_0.7\_9 &          97 &      86 &          79 &       23.71 &   21.21 &       27.41 &           94.59 &  94.92 &      94.59 \\
 BPP\_200\_100\_0.1\_0.8\_1 &         117 &     112 &         112 &       48.07 &   46.64 &       69.57 &           92.69 &  92.90 &      92.97 \\
 BPP\_200\_100\_0.1\_0.7\_1 &          69 &      65 &          61 &       20.25 &   19.44 &       27.35 &           78.83 &  78.92 &      79.06 \\
  BPP\_200\_75\_0.1\_0.7\_4 &          66 &      48 &          67 &       12.80 &    8.78 &       19.97 &           80.44 &  81.72 &      80.37 \\
  BPP\_200\_75\_0.2\_0.7\_1 &          59 &      43 &          46 &        8.50 &    5.70 &        9.00 &           86.87 &  87.82 &      86.80 \\
 BPP\_750\_300\_0.1\_0.7\_7 &         230 &     215 &         231 &      465.58 &  431.99 &      761.27 &          299.97 & 299.85 &     298.28 \\
 BPP\_750\_300\_0.1\_0.7\_8 &         288 &     215 &         252 &      642.50 &  438.36 &      868.59 &          301.65 & 302.74 &     301.51 \\
 BPP\_750\_300\_0.1\_0.8\_2 &         385 &     329 &         331 &     1222.12 &  984.78 &     1575.72 &          344.83 & 346.43 &     344.20 \\
 BPP\_750\_300\_0.1\_0.8\_3 &         330 &     327 &         316 &      987.16 &  992.71 &     1487.61 &          333.55 & 334.33 &     333.36 \\
 BPP\_750\_300\_0.1\_0.8\_7 &         397 &     317 &         321 &     1238.79 &  900.67 &     1466.03 &          334.83 & 336.11 &     335.26 \\
 BPP\_750\_300\_0.2\_0.7\_6 &         377 &     194 &         226 &      924.06 &  348.74 &      695.31 &          339.59 & 341.83 &     343.53 \\
 BPP\_750\_300\_0.2\_0.7\_7 &         302 &     202 &         241 &      652.23 &  374.72 &      761.56 &          346.00 & 346.91 &     345.16 \\
 BPP\_750\_300\_0.2\_0.7\_8 &         266 &     190 &         213 &      555.75 &  343.14 &      643.76 &          339.14 & 340.80 &     340.28 \\
 BPP\_750\_300\_0.1\_0.7\_5 &         262 &     225 &         243 &      676.38 &  553.77 &      991.46 &          302.62 & 305.25 &     302.32 \\
 BPP\_750\_300\_0.1\_0.7\_2 &         282 &     218 &         225 &      607.61 &  424.71 &      723.86 &          308.22 & 310.59 &     307.71 \\
 BPP\_750\_300\_0.2\_0.7\_1 &         251 &     215 &         213 &      393.61 &  324.76 &      513.24 &          343.90 & 343.73 &     344.31 \\
 BPP\_750\_300\_0.1\_0.7\_4 &         252 &     213 &         248 &      526.86 &  424.10 &      839.43 &          299.13 & 301.11 &     298.10 \\
 BPP\_750\_300\_0.1\_0.7\_9 &         251 &     266 &         255 &      518.30 &  576.85 &      867.30 &          303.05 & 301.99 &     301.48 \\
 BPP\_750\_300\_0.2\_0.7\_0 &         264 &     199 &         240 &      436.00 &  291.99 &      602.91 &          338.58 & 337.66 &     336.33 \\
 BPP\_750\_300\_0.2\_0.7\_2 &         330 &     206 &         284 &      595.43 &  289.76 &      719.85 &          344.99 & 347.38 &     345.10 \\
 BPP\_750\_300\_0.2\_0.7\_3 &         291 &     223 &         258 &      496.82 &  347.74 &      667.53 &          333.32 & 332.97 &     333.40 \\
 BPP\_750\_300\_0.2\_0.7\_4 &         304 &     201 &         224 &      539.88 &  299.66 &      567.41 &          343.59 & 343.89 &     343.54 \\
 BPP\_750\_300\_0.1\_0.7\_6 &         254 &     221 &         224 &      545.02 &  452.20 &      755.89 &          302.22 & 303.04 &     300.90 \\
 BPP\_750\_300\_0.2\_0.7\_0 &         264 &     199 &         240 &      437.81 &  294.21 &      618.61 &          338.58 & 337.66 &     336.33 \\
 BPP\_750\_300\_0.2\_0.7\_2 &         330 &     206 &         284 &      595.81 &  298.91 &      749.31 &          344.99 & 347.38 &     345.10 \\
 BPP\_750\_300\_0.2\_0.7\_5 &         239 &     190 &         209 &      398.43 &  271.92 &      508.35 &          335.84 & 337.90 &     335.33 \\
\bottomrule
\end{tabular}}
\end{table*}

\end{document}